\documentclass[12pt]{amsart} 
\usepackage{amsmath, amsfonts, amssymb, latexsym, amscd, enumerate}
\usepackage{pictexwd,dcpic}
\usepackage{graphicx}

\newtheorem{theorem}{Theorem}[section]
\newtheorem{lemma}[theorem]{Lemma}
\newtheorem{proposition}[theorem]{Proposition}
\newtheorem{corollary}[theorem]{Corollary}
\newtheorem{remark}[theorem]{Remark}
\newtheorem{example}[theorem]{Example}
\newtheorem{definition}[theorem]{Definition}
\newtheorem{notation}[theorem]{Notation}

\numberwithin{equation}{section}

\def\zN{\mathbb N}
\def\zZ{\mathbb Z}
\def\zQ{\mathbb Q}
\def\zR{\mathbb R}
\def\zC{\mathbb C}
\def\zT{\mathbb T}

\def\cE{\mathcal E}
\def\cF{\mathcal F}
\def\cL{\mathcal L}
\def\cK{\mathcal K}
\def\cI{\mathcal I}
\def\cO{\mathcal O}
\def\cT{\mathcal T}

\def\fI{\mathfrak I}
\def\sq{$\text{}\hfill\square$}

\def\pf{\it Proof. \rm}
\def\st{\,\bigl\vert\,}

\newcommand{\ds}[1]{\displaystyle#1}
\newcommand{\abs}[1]{\vert#1\vert}
\newcommand{\norm}[1]{\vert\vert#1\vert\vert}
\newcommand{\innprod}[2]{\langle#1, #2\rangle}
\newcommand{\mat}[4]{\begin{pmatrix}#1&#2\\#3&#4\end{pmatrix}}
\newcommand{\diag}[1]{\mbox{Diag}(#1)}

\title[$C^*$-algebras associated with topological group quivers II]{$C^*$-algebras associated with topological group quivers II:\\ \rm $K$-groups}
\author{Shawn J. $\rm M^\MakeLowercase{c}$Cann}
\begin{document}
\maketitle

\begin{abstract} Topological quivers generalize the notion of directed graphs in which
the sets of vertices and edges are locally compact (second countable) Hausdorff spaces. 
Associated to a topological quiver $Q$ is a $C^*$-correspondence, and in turn, a Cuntz-Pimsner 
algebra $C^*(Q).$ Given $\Gamma$ a locally compact group and $\alpha$ and $\beta$ 
endomorphisms on $\Gamma,$ one may construct a topological quiver $Q_{\alpha,\beta}(\Gamma)$
with vertex set $\Gamma,$ and edge set $\Omega_{\alpha,\beta}(\Gamma)=
\{(x,y)\in\Gamma\times\Gamma\st \alpha(y)=\beta(x)\}.$ In \cite{Mc1}, the author
examined the Cuntz-Pimsner algebra $\cO_{\alpha,\beta}(\Gamma):=C^*(Q_{\alpha,\beta}(\Gamma))$
and found  generators (and their relations) of $\cO_{\alpha,\beta}(\Gamma).$ In this paper, the author uses
this information to create a six term exact sequence in order to calculate the $K$-groups of $\cO_{\alpha,\beta}(\Gamma).$
\end{abstract}

\section{Introduction and Notation}
	\subsection{Background}Given a quintuple $Q = (X,E, r, s, \lambda),$ where $X$ and $E$ are locally compact (second countable) Hausdorff
spaces, $r$ and $s$ are continuous maps from $X$ to $E$ with $r$ open, and $\lambda = \{\lambda_x\}_{x\in E}$
is a system of Radon measures, one can create a corresponding Cuntz-Pimsner $C^*$-algebra $C^*(Q).$
In \cite{EaHR}, Exel, an Huef and Raeburn define $C^*$-algebras associated with a system $(B,\alpha,L)$
where $\alpha$ is an endomorphism of a unital $C^*$-algebra $B$ and $L$ is a positive linear map $L:B\to B$ such
that $L(\alpha(a)b) = aL(b)$ for all $a, b\in B$ called a \emph{transfer operator}. In fact, the $C^*$-algebra they generate is a
Cuntz-Pimsner algebra and under certain restrictions, a $C^*$-algebra associated with a topological quiver;
in particular, when $B=C(\zT^d)$ the continuous function on the $d$-torus,
$F\in M_d(\zZ)$ and $\alpha$ is the endomorphism 
$$\alpha(f)(e^{2\pi it})=f(e^{2\pi i Ft})$$
for $f\in C(\zT^d)$ and $t\in\zR^d.$ 
Exel, an Huef and Raeburn then determine a six term exact sequence in which to use to calculate the $K$-groups
of these $C^*$-algebras. In \cite{Mc1}, the author considers a certain class of topological quivers (which extend the notions
of Exel, an Huef and Raeburn) $Q=(\Gamma,\Omega_{\alpha,\beta}(\Gamma),r,s,\lambda)$ where $\Gamma$ is a locally
compact group, $\alpha$ and $\beta$ are endormorphism of $\Gamma,$
$$\Omega_{\alpha,\beta}(\Gamma)=\{(x,y)\in\Gamma\times\Gamma\st \alpha(y)=\beta(x)\}$$
and $\lambda$ is an appropriate family of Radon measures. The resulting Cuntz-Pimsner $C^*$-algebra, 
denoted $\cO_{\alpha,\beta}(\Gamma),$ was then examined and certain generators and relations where found.
We now proceed to generalize the six term exact sequence considered in \cite{EaHR} to $C^*$-algebras of the form
$\cO_{\alpha,\beta}(\Gamma)$ where $\Gamma$ is a compact group.

	\subsection{Notation}The sets of natural numbers, integers, rationals numbers, real numbers and complex numbers will be denoted by
$\zN$, $\zZ$, $\zQ$, $\zR$, and $\zC,$ respectively. Convention: $\zN$ does not contain zero.
$\zZ_0^+$ will denote the set $\zN\cup\{0\},$ $\zR^+$ denotes the set $\{r\in\zR\st r>0\}$ and $\zR_0^+=\zR^+\cup\{0\}.$
Finally, $\zZ_p$ denotes the abelian group $\zZ/p\zZ=\{0,1,...,p-1\mod p\}$ and 
$\zT$ denotes the torus $\{z\in\zC\st \abs{z}=1\}.$ Whenever convenient, view $\zZ_p\subset\zT$ by 
$\zZ_p\cong\{z\in\zT\st z^p=1\}.$

For a topological space $Y$, the closure of $Y$ is denoted $\overline{Y}.$ Given a locally compact Hausdorff space $X$, let
\begin{enumerate} 
\item $C(X)$ be the continuous complex functions on $X$;\index{C$(X)$}
\item $C_b(X)$ be the continuous and bounded complex functions on $X$;\index{C$\mbox{}_b(X)$}
\item $C_0(X)$ be the continuous complex functions on $X$ vanishing at infinity;\index{C$\mbox{}_0(X)$}
\item $C_c(X)$ be the continuous complex functions on $X$ with compact support.\index{C$\mbox{}_c(X)$}
\end{enumerate}
The supremum norm is denoted $\norm{\cdot}_\infty$ and defined by
$$\norm{f}_\infty=\sup_{x\in X}\{\abs{f(x)}\}$$
for each continuous map $f:X\to\zC.$ For a continuous function $f\in C_c(X),$ denote the open support of $f$ by 
$\mbox{osupp }f=\{x\in X\st f(x)\ne 0\}$ and the support of $f$ by $\mbox{supp }f=\overline{\mbox{osupp} f}.$  

For $C^*$-algebras $A$ and $B$, $A$ is isomorphic to $B$ will be written $A\cong B;$ for example, we use
$C(\zT^d)\otimes M_{N}(\zC)\cong M_{N}(C(\zT^d)).$ 
Moreover, $A^{\oplus n}$ denotes the $n$-fold direct sum $A\oplus\cdots\oplus A.$ 
Given a group $\Gamma$ and a ring $R$, a normal subgroup, $N$, of $\Gamma$ is denoted $N\lhd\Gamma$ and
an ideal, $I$, of $R$ is denoted $I\lhd R.$ Note if $R$ is a $C^*$-algebra then the term ideal denotes a closed two-sided ideal. 
Furthermore, $\mbox{End}(\Gamma)$ ($\mbox{End}(R)$) and $\mbox{Aut}(\Gamma)$
($\mbox{Aut}(R)$) denotes the set of endomorphisms of $\Gamma$ ($R$) and automorphisms of $\Gamma$ $(R$), respectively. 
For a map $\gamma:\Gamma\to\mbox{Aut}(A),$ the \emph{fixed point set}\index{Fixed Point Set} 
is denoted $A^\gamma$ and defined by 
$$A^\gamma=\{a\in A\st \gamma(g)(a)=a\mbox{ for each }g\in \Gamma\}.$$

Let $\alpha\in C(X)$ then $\alpha^\#\in\mbox{End(C(X))}$ denotes the endomorphism of $C(X)$ defined by
$$\alpha^\#(f)=f\circ\alpha\qquad\mbox{for each $f\in C(X)$}.$$
Let $S$ be a set and define the Kronecker delta function $\delta:S\times S\to\{0,1\}$ by
$$\delta_s^r:=\delta(s,r)=\begin{cases}
0&\mbox{if }s\ne r\\
1&\mbox{if }s=r
\end{cases}.$$

The set of $n$ by $n$ matrices with coefficients in a set $R$ will be denoted $M_n(R)$ and for any $F\in M_n(R),$
the transpose of $F$ is denoted $F^T$. Given a function $\sigma: R\to S$, we may create an augmented function
$\sigma_n:M_n(R)\to M_n(S)$ via
$$\sigma_n((r_{i,j})_{i,j=1}^n)=(\sigma(r_{i,j}))_{i,j=1}^n$$
for each $(r_{i,j})_{i,j=1}^n\in M_n(R).$  Given vectors $v=(v_1,...,v_n)$ of length $n$ and $w=(w_1,...,w_m)$ of length $m$, denote $(v,w)$ to be the vector $(v,w)=(v_1,...,v_n,w_1,...,w_m)$ of length $n+m.$
\section{Preliminairies}\subsection{Hilbert $C^*$-modules}We begin by defining Hilbert $C^*$-modules. Further details and references can be found in \cite{lan, RW}.

\begin{definition}\label{Hbmod}\cite{lan} \rm  If $A$ is a $C^*$-algebra, then a \emph{(right) Hilbert $A$-module}\index{Hilbert $C^*$-module} is a Banach space $\cE_A$
together with a right action of $A$ on $\cE_A$ and an $A$-valued inner product $\innprod{\cdot}{\cdot}_A$
satisfying
\begin{enumerate}
\item $\innprod{\xi}{\eta a}_A=\innprod{\xi}{\eta}_A a$
\item $\innprod{\xi}{\eta}_A =\innprod{\eta}{\xi}_A^*$
\item $\innprod{\xi}{\xi}\ge 0$ and $\norm{\xi}=\norm{\innprod{\xi}{\xi}_A^{1/2}}_A$
\end{enumerate}
for all $\xi$, $\eta\in\cE_A$ and $a\in A$ (if the context is clear, we denote $\cE_A$ simply by $\cE$). 
For Hilbert $A$-modules $\cE$ and $\cF$, call a function $T:\cE\to\cF$ \emph{adjointable}
\index{Hilbert $C^*$-module! Adjointable Operator, $\cL(\cE,\cF)$}
\index{Adjointable Operator, $\cL(\cE,\cF)$}
if there is a function $T^*:\cF\to\cE$ such that
$\innprod{T(\xi)}{\eta}_A=\innprod{\xi}{T^*(\eta)}_A$ for all $\xi\in\cE$ and $\eta\in\cF$.
Let $\cL(\cE,\cF)$ denote the set of adjointable ($A$-linear) operators from $\cE$ to $\cF$. 
If $\cE=\cF$, then $\cL(\cE):=\cL(\cE,\cE)$ is a $C^*$-algebra (see \cite{lan}.)
Let $\cK(\cE, \cF)$ denote the closed two-sided ideal of \emph{compact operators}
\index{Hilbert $C^*$-module! Compact Operators, $\cK(\cE,\cF)$} 
\index{Compact Operators, $\cK(\cE,\cF)$}
given by
$$\cK(\cE,\cF):=\overline{\mbox{span}}\{\theta_{\xi,\eta}^{\cE,\cF}\st\xi\in\cE,\,\eta\in\cF\}$$
where $$\theta_{\xi,\eta}^{\cE,\cF}(\zeta)=\xi\innprod{\eta}{\zeta}_A\qquad\mbox{for each $\zeta\in\cE$}.$$ 
Similarly, $\cK(\cE):=\cK(\cE,\cE)$ and $\theta_{\xi,\eta}^\cE$ (or $\theta_{\xi,\eta}$ if understood) denotes 
$\theta_{\xi,\eta}^{\cE,\cE}$.
For Hilbert $A$-module $\cE$, the linear span of $\{\innprod{\xi}{\eta}\st\xi,\eta\in\cE\}$, denoted $\innprod{\cE}{\cE}$,
once closed is a two-sided ideal of $A$. Note that $\cE\innprod{\cE}{\cE}$ is dense in $\cE$ (\cite{lan}). 
The Hilbert module $\cE$ is called \emph{full}\index{Hilbert $C^*$-module! Full} if $\innprod{\cE}{\cE}$ is dense in $A$. 
The Hilbert module $A_A$ refers to the Hilbert module $A$ over itself, where $\innprod{a}{b}=a^*b$ for all $a,b\in A$.

An \emph{algebraic generating set}\index{Hilbert $C^*$-module! Algebraic Generating Set} for $\cE$ is a subset $\{u_i\}_{i\in\cI}\subset \cE$ for some indexing set $\cI$ such that $\cE$ equals the linear span of $\{u_i\cdot a\st\ i\in\cI, a\in A\}.$ 
\end{definition}

\begin{definition}\label{ONB}\cite{KW} \rm A subset $\{u_i\}_{i\in\cI}\subset \cE$ is called a \emph{basis}\index{Basis}
provided the following reconstruction formula holds for all $\xi\in\cE:$
$$ \xi=\sum_{i\in\cI} u_i\cdot\innprod{u_i}{\xi}\qquad(\mbox{in }\cE,\norm{\cdot}.)$$
If $\innprod{u_i}{u_j}=\delta_i^j$ as well, call $\{u_i\}_{i\in\cI}$ an \emph{orthonormal basis}
\index{Basis! Orthonormal} of $\cE$.
\end{definition}  

\begin{remark}\rm The preceding definition is in accordance with the finite version in \cite{KW}, 
but many other versions exist such as in
\cite{EaHR} where $\{u_i\}_{i=1}^n$ is called a finite Parseval frame, or in \cite{Yam} where this is taken as the definition
for \it finitely generated.\rm\index{Hilbert $C^*$-module! Finitely Generated} 
There has been substantial work done on similar frames (see \cite{HJLM}).
\end{remark}

The following notions of $C^*$-correspondence and morphism may be found in \cite{MT, BB3, BB4, BB5, EaHR, FLR, FMR, K1}

\begin{definition}\label{Hbmorph}\cite{BB4, BB5} \rm If $A$ and $B$ are $C^*$-algebras, then an \emph{$A-B$ $C^*$-correspondence}\index{C$\mbox{}^*$-correspondence} $\cE$  is a right Hilbert 
$B$-module $\cE_B$ together with a left action of $A$ on $\cE$ given by a $*$-homomorphism $\phi_A:A\to\cL(\cE)$, 
$a\cdot\xi=\phi_A(a)\xi$ for $a\in A$ and $\xi\in\cE$. We may occasionally write, $_A\cE_B$ to denote an $A-B$ 
$C^*$-correspondence and $\phi$ instead of $\phi_A$.
Furthermore, if $_{A_1}\cE_{B_1}$ and $_{A_2}\cF_{B_2}$ are $C^*$-correspondences, then
a \emph{morphism}\index{C$\mbox{}^*$-correspondence! Morphism} 
$(\pi_1, T, \pi_2):\cE\to\cF$ consists of $*$-homomorphisms $\pi_i:A_i\to B_i$ and a linear map
$T:\cE\to\cF$ satisfying
\begin{enumerate}
\item[(i)] $\pi_2(\innprod{\xi}{\eta}_{A_2})=\innprod{T(\xi)}{T(\eta)}_{B_2}$
\item[(ii)] $T(\phi_{A_1}(a_1)\xi)=\phi_{B_1}(\pi_1(a_1))T(\xi)$
\item[(iii)] $T(\xi)\pi_2(a_2)=T(\xi a_2)$
\end{enumerate} for all $\xi,\eta\in\cE$ and $a_i\in A_i$.
\end{definition}

\begin{notation}\rm When $A=B$, we refer to $_A\cE_A$ as a $C^*$-correspondence over $A$. For $\cE$ a 
$C^*$-correspondence over $A$ and $\cF$ a $C^*$-correspondence over $B$, a morphism 
 $(\pi,T,\pi):\cE\to\cF$ will be denoted by $(T,\pi)$. 
\end{notation}

\begin{definition}\cite{MT} \rm If $\cF$ is the Hilbert module $ _CC_C$ where $C$ is a $C^*$-algebra with the inner product 
$\innprod{x}{y}_B=x^*y$ then call a morphism $(T,\pi):$ $_A\cE_B\to C$ of Hilbert
modules a \emph{representation}\index{C$\mbox{}^*$-correspondence! Representation} of $_A\cE_B$ into $C.$ 
\end{definition}

\begin{remark}\rm Note that a representation of $_A\cE_B$ need only satisfying $(i)$ and $(ii)$ 
of definition \ref{Hbmorph} as it was unnecessary to require (iii). For a proof, see \cite[Remark 2.7]{Mc1}.
\end{remark}

A morphism of Hilbert modules $(T,\pi):\cE\to\cF$ yields a $*$-homomorphism $\Psi_T:\cK(\cE)\to\cK(\cF)$ by
$$\Psi_T(\theta_{\xi,\eta}^\cE)=\theta_{T(\xi),T(\eta)}^\cF$$
for $\xi,\eta\in\cE$ and if $(S,\sigma):\mathcal D\to\cE$, and $(T,\pi):\cE\to\cF$ are morphisms of Hilbert modules then
$\Psi_T\circ\Psi_S=\Psi_{T\circ S}$. In the case where $\cF=B$ a $C^*$-algebra, we may first identify $\cK(B)$ as $B$,
and a representation $(T,\pi)$ of $\cE$ in a $C^*$-algebra $B$ yields a $*$-homomorphism $\Psi_T:\cK(\cE)\to B$ given
by $$\Psi_T(\theta_{\xi,\eta})=T(\xi)T(\eta)^*.$$

\begin{definition}\cite{MT} \rm For a $C^*$-correspondence $\cE$ over $A$, denote the ideal $\phi^{-1}(\cK(\cE))$ of $A$ by $J(\cE),$\index{J$(\cE)$} and let $J_\cE=J(\cE)\cap(\ker \phi)^\perp$\index{J$\mbox{}_\cE$} where $(\ker\phi)^\perp $ is the 
ideal $\{a\in A\st ab=0 \mbox{ for all }b\in\ker\phi\}$ .
If $_A\cE_A$ and $_B\cF_B$ are $C^*$-correspondences over $A$ and $B$ respectively and $K\lhd J(\cE)$, a morphism
$(T,\pi):\cE\to\cF$ is called \emph{coisometric on $K$}\index{C$\mbox{}^*$-correspondence! Representation! Coisometric on $K$} if $$\Psi_T(\phi_A(a))=\phi_B(\pi(a))$$
for all $a\in K$, or just \emph{coisometric},\index{C$\mbox{}^*$-correspondence! Representation! Coisometric} if $K=J(\cE)$.
\end{definition}

\begin{notation}\rm We denote $C^*(T,\pi)$ to be the $C^*$-algebra generated by $T(\cE)$ and $\pi(A)$ where
$(T,\pi):\cE\to B$ is a representation of $_A\cE_A$ in a $C^*$-algebra $B$. Furthermore, if $\rho:B\to C$ is a 
$*$-homomorphism of $C^*$-algebras, then $\rho\circ (T,\pi)$ denotes the representation $(\rho\circ T,\rho\circ\pi)$ of $\cE$.
\end{notation}

\begin{definition}\label{CPDefs}\cite{MT} \rm A morphism $(T_\cE,\pi_\cE)$ coisometric on an ideal $K$ is said to be \emph{universal}\index{C$\mbox{}^*$-correspondence! Representation! Universal} 
if whenever $(T,\pi):\cE\to B$ is a representation coisometric on $K$, there exists a $*$-homomorphism 
$\rho:C^*(T_\cE,\pi_\cE)\to B$ with $(T,\pi)=\rho\circ(T_\cE,\pi_\cE)$.  The universal $C^*$-algebra 
$C^*(T_\cE,\pi_\cE)$ is called the \emph{relative Cuntz-Pimsner algebra}\index{Cuntz-Pimsner Algebra}\index{$\cO(K,\cE)$}
\index{$\cO_\cE$} 
of $\cE$ determined by the ideal $K$ and 
denoted by $\cO(K,\cE)$. If $K=0$, then $\cO(K,\cE)$ is denoted
by $\cT(\cE)$ and called the \emph{universal Toeplitz $C^*$-algebra}\index{Toeplitz-Pimsner Algebra}\index{$\cT(\cE)$} for $\cE$. We denote $\cO(J_\cE,\cE)$ by $\cO_\cE$. 
\end{definition}
	\subsection{Topological Quivers}\begin{definition}\label{TQ}\cite{MT} \rm A \emph{topological quiver}\index{Topological Quiver} 
(or \emph{topological directed graph}\index{Topological Directed Graph}) $Q=(X,E,Y,r,s,\lambda)$ is a diagram
$$\begindc{\commdiag}[5]
\obj(10,0){$E$}
\obj(0,0){$X$}
\obj(20,0){$Y$}
\mor{$E$}{$X$}{$s$}[-1,0]
\mor{$E$}{$Y$}{$r$}
\enddc$$
where $X,E,$ and $Y$ are second countable locally compact Hausdorff spaces, $r$ and $s$ are continuous maps with $r$ open, along with a family $\lambda=\{\lambda_y\vert y\in Y\}$ of Radon measures on $E$ satisfying
\begin{enumerate}
\item $\mbox{supp }\lambda_y=r^{-1}(y)$ for all $y\in Y$, and
\item $y\mapsto\lambda_y(f)=\int_Ef(\alpha)d\lambda_y(\alpha)\in C_c(Y)$ for $f\in C_c(E).$
\end{enumerate}
\end{definition}

\begin{remark}\rm If $X=Y$ then write $Q=(X,E,r,s,\lambda)$ in lieu of $(X,E,X,r,s,\lambda).$
\end{remark}

\begin{remark}\rm The author provides a broad history and a series of examples of topological quivers in \cite{McThesis, Mc1}.
\end{remark}

Given a topological quiver $Q=(X,E,Y,r,s,\lambda)$, one may associate a correspondence  
$\cE_Q$ of the $C^*$-algebra $C_0(X)$ to the $C^*$-algebra $C_0(Y)$. Define left and right actions 
$$(a\cdot\xi\cdot b)(e)=a(s(e))\xi(e)b(r(e))$$
by $C_0(X)$ and $C_0(Y)$ respectively on $C_c(E)$. Furthermore, define the $C_c(Y)$-valued inner product
$$\innprod{\xi}{\eta}(y)=\int_{r^{-1}(y)}\overline{\xi(\alpha)}\eta(\alpha)d\lambda_y(\alpha)$$ 
for $\xi,\eta\in C_c(E)$, $y\in Y,$ and let $\cE_Q$\index{C$\mbox{}^*$-correspondence! Associated with a Topological Quiver} be the completion of $C_c(E)$ with respect to the norm 
$$\norm{\xi}=\norm{\innprod{\xi}{\xi}^{1/2}}_\infty=\norm{\lambda_y(\abs{\xi}^2)}_\infty^{1/2}.$$ 

\begin{definition}\rm Given topological quiver $Q$ over a space $X$, define the $C^*$-algebra, $C^*(Q)$
\index{C$\mbox{}^*$-algebra Associated with a Topological Quiver}\index{C$\mbox{}^*(Q)$}
\index{Topological Quiver! $C^*$-algebra Associated with} 
associated with $Q$ to be the Cuntz-Pimnser $C^*$-algebra $\cO_{\cE_Q}$\index{Cuntz-Pimsner Algebra} 
of the correspondence $\cE_Q$ over $A=C_0(X)$.
\end{definition}

	\subsection{Topological Group Quivers}\begin{definition}\label{TopGrpQuiver}\rm Let $\Gamma$ be a (second countable) locally compact group and let $\alpha,\beta\in\mbox{End}(\Gamma)$ be continuous. Define the closed subgroup, $\Omega_{\alpha,\beta}(\Gamma),$ 
of $\Gamma\times\Gamma,$\index{$\Omega_{\alpha,\beta}(\Gamma)$}
$$\Omega_{\alpha,\beta}(\Gamma)=\{(x,y)\in \Gamma\times\Gamma\st \alpha(y)=\beta(x)\}$$
and let $Q_{\alpha,\beta}(\Gamma)=(\Gamma,\Omega_{\alpha,\beta}(\Gamma), r,s,\lambda)$
\index{Q$\mbox{}_{\alpha,\beta}(\Gamma)$}\index{Topological Quiver}\index{Topological Quiver! $C^*$-algebra Associated with}  where $r$ and $s$ are the group homomorphisms defined by 
$$r(x,y)=x \qquad\mbox{and}\qquad s(x,y)=y$$ 
for each $(x,y)\in\Omega_{\alpha,\beta}(\Gamma)$ and $\lambda_x$ for $x\in\Gamma$ is the measure on 
$$r^{-1}(x)=\{x\}\times\alpha^{-1}(\beta(x))$$ 
defined by 
$$\lambda_x(B)=\mu(y^{-1}s(B\cap r^{-1}(x))\cap\ker\alpha)\qquad\mbox{(for any}\, y\in\alpha^{-1}(\beta(x)))$$
for each measurable $B\subseteq\Omega_{\alpha,\beta}(\Gamma)$ 
where $\mu$ is a left Haar measure (normalized if possible) on $r^{-1}(1_\Gamma)=\{1\}\times\ker\alpha$ 
(a closed normal subgroup of $\Gamma\times\Gamma;$ hence, a locally compact group). Note if $r^{-1}(x)=\emptyset$
then $\alpha^{-1}(\beta(x))=\emptyset$ and so $\lambda_x=0.$ 
This measure is well-defined, 
$$\mbox{supp }\lambda_x=\{x\}\times y\ker\alpha=\{x\}\times \alpha^{-1}(\beta(x))=r^{-1}(x)$$
and $y\mapsto \lambda_y(f)$ is a continuous compactly supported function (cf. \cite[Definition 3.1]{Mc1}.

Call $Q_{\alpha,\beta}(\Gamma)$ a \emph{topological group relation.}\index{Topological Group Relation}\index{Topological Relations}
Define $\cE_{\alpha,\beta}(\Gamma)$\index{$\cE_{\alpha,\beta}(\Gamma)$}\index{C$\mbox{}^*$-correspondence} 
to be the $C_0(\Gamma)$-correspondence $\cE_{Q_{\alpha,\beta}(\Gamma)}$ and
form the Cuntz-Pimsner algebra\index{$\cO_{\alpha,\beta}(\Gamma)$! Definition}\index{Cuntz-Pimsner Algebra}
$$\cO_{\alpha,\beta}(\Gamma):=C^*(Q_{\alpha,\beta}(\Gamma))=\cO(J_{\cE_{\alpha,\beta}(\Gamma)},\cE_{\alpha,\beta}(\Gamma))$$
and the Toeplitz-Pimsner algebra\index{$\cT_{\alpha,\beta}(\Gamma)$! Definition}\index{Toeplitz-Pimsner Algebra}
$$\cT_{\alpha,\beta}(\Gamma):=\cT(Q_{\alpha,\beta}(\Gamma)).$$
\end{definition}

\begin{remark}\rm It will be implicitly assumed that $\Gamma$ is second countable. Furthermore, since $\Gamma$ is locally compact Hausdorff, $r^{-1}(x)$ is closed and locally compact. Moreover, whenever $r$ is a local homeomorphism, $r^{-1}(x)$ is discrete and hence, $\lambda_x$ is counting measure (normalized when $\abs{\ker\alpha}<\infty$.)
\end{remark}\begin{example}[\cite{Mc1}]\label{ToriQuiver}\rm\index{Topological Group Relation! $\zT^d$-quivers}\index{Topological Group Relation}
\index{Topological Relations}\index{Topological Quiver}\index{$\cO_{F,G}(\zT^d)$! Definition}
For the compact abelian group $\zT^d,$ note $\mbox{End}(\zT^d)\cong M_d(\zZ)$ (\cite{W}); that is, 
an element $\sigma\in\mbox{End}(\zT^d)$ is of the form $\sigma_F$ for some $F\in M_d(\zZ)$ where 
$$\sigma_F(e^{2\pi it})=e^{2\pi i Ft}$$\index{$\sigma_F$}
for each $t\in\zZ^d.$ To simplify notation, use $F$ and $G$ in place of $\sigma_F$ and $\sigma_G$ whenever convenient.
For instance,
$$Q_{F,G}(\zT^d):=Q_{\sigma_{F},\sigma_{G}}(\zT^d)$$
and the $C^*$-correspondence
$$\cE_{F,G}(\zT^d):=\cE_{\sigma_F,\sigma_G}(\zT^d)$$
where $F,G\in M_d(\zZ)$. We will consider the cases when these maps are surjective; that is,
$\det F$ and $\det G$ are non-zero.

Let $F, G\in M_d(\zZ)$ where $\det F,\det G\ne 0$. Then $\abs{\ker\sigma_F}=\abs{\det F}$
and so, the $C(\zT^d)$-valued inner product becomes
$$\innprod{\xi}{\eta}(x)=\frac{1}{\abs{\det F}}\sum_{\sigma_{F}(y)=\sigma_{G}(x)}\overline{\xi(x,y)}\eta(x,y)$$
for $\xi,\eta\in\cE_{F,G}(\zT^d)$ and $x\in\zT^d.$
This is a finite sum since the number of solutions, $y,$ 
to $\sigma_F(y)=\sigma_G(x)$ given any $x\in\zT^d$ is $\abs{\det F}<\infty.$
\end{example}

\begin{remark}\rm The left action, $\phi,$ is defined by
$$\phi(a)\xi(x,y)=a(y)\xi(x,y)$$
for $a\in C(\zT^d),$ $\xi\in C(\Omega_{F,G}(\zT^d))$ and $(x,y)\in\Omega_{F,G}(\zT^d).$  
Note: $\phi$ is injective. To see this claim, let $a\in C(\zT^d)$ and assume $\phi(a)\xi=0$ for each $\xi\in C(\Omega_{F,G}(\zT^d)).$ Then $a(y)\xi(x,y)=0$ for each $(x,y)\in\Omega_{F,G}(\zT^d)$ and $\xi\in C(\Omega_{F,G}(\zT^d)).$ Since 
$s(\Omega_{F,G}(\zT^d))=\{y\in\zT^d\st (x,y)\in\Omega_{F,G}(\zT^d)\}=\zT^d$ by the surjectivity of $\sigma_F,$
$a=0.$ 
\end{remark}

\begin{remark}\rm It was shown in \cite{Mc1} that one may assume the matrix $F$ is positive diagonal. 
\end{remark}

Let $F=\diag{a_1,...,a_d}\in M_d(\zZ), G=(b_{jk})_{j,k=1}^d\in M_d(\zZ)$ where $a_j>0$ for each $j=1,...,d,$ $\det G\ne0$ and
let $G_j$ denote the $j$-th row of $G$, $(b_{jk})_{k=1}^d.$ Further, let $N=\det F=\prod_{j=1}^d a_j>0$ and let $$\fI(F)=\{\nu=(\nu_j)_{j=1}^d\in\zZ^d\st 0\le \nu_j\le a_j-1\}.\index{$\fI(F)$}$$ 
The $C(\zT^d)$-valued  inner product becomes
$$\innprod{\xi}{\eta}(x)=\frac{1}{N}\sum_{\sigma_F(y)=\sigma_G(x)} \overline{\xi(x,y)}\eta(x,y)$$
for all $\xi,\eta\in C(\Omega_{F,G}(\zT^d))$ and $x\in \zT^d.$

Given $\nu\in\fI(F)$, define $u_\nu\in C(\Omega_{F,G}(\zT^d))$ by 
$$u_\nu(x,y)=y^\nu=\prod_{j=1}^d y^{\nu_j}$$
for $(x,y)\in \Omega_{F,G}(\zT^d).$ It was shown in \cite{Mc1} that $\{u_\nu\}_{\nu\in\fI(F)}$ is a basis for $\cE_{F,G}(\zT^d)$
and also the following:

\begin{theorem}\cite[Theorem 3.23]{Mc1}\label{TdQuiverGen}\rm \index{$\cO_{F,G}(\zT^d)$}\index{$\cO_{\alpha,\beta}(\Gamma)$! Presentation}
\index{Cuntz-Pimsner Algebra}
Let $F=\diag{a_1,...,a_d}, G\in M_d(\zZ)$ where $\det F, \det G\ne0$ and let $G_j$ be the $j$-th row vector of $G$. Further, let $\fI(F)$ denote the set $\{\nu=(\nu_j)_{j=1}^d\in\zZ^d\st 0\le \nu_j\le a_j-1\}$. Then
$\cO_{F,G}(\zT^d)$ is the universal $C^*$-algebra generated by isometries $\{S_\nu\}_{\nu\in\fI(F)}$ and (full spectrum) commuting unitaries $\{U_j\}_{j=1}^d$ that satisfy the relations
\begin{enumerate}
\item $S_\nu^*S_{\nu^\prime}=\innprod{u_\nu}{u_{\nu^\prime}}=\delta_{\nu}^{\nu^\prime},$
\item $U^\nu S=S_\nu$ for all $\nu\in\fI(F),$
\item $U_j^{a_j}S=SU^{G_j},$ for all $j=1,...,d$ and
\item $1=\sum_{\nu\in\fI(F)} S_\nu S_\nu^*=\sum_{\nu\in\fI(F)} U^\nu SS^*U^{-\nu}$
\end{enumerate}
where $U^\nu$ denotes $\prod_{j=1}^dU_j^{\nu_j}.$ Furthermore, $\cT_{\alpha,\beta}(\Gamma)$ is the universal $C^*$-algebra generated by isometries $\{S_\nu\}_{\nu\in\fI(F)}$ and commuting unitaries $\{U_j\}_{j=1}^d$ that satisfy relations (1)-(3)
\end{theorem}
\section{Six Term Exact Sequence for $\cO_{\alpha,\beta}(\Gamma)$}In this section, we follow and extend the approach of \cite{EaHR} to create a six term exact sequence.
Let $\Gamma$ be a compact group with $\alpha,\beta\in\mbox{End}(\Gamma).$
Suppose the left action for the correspondence, $\phi,$ is injective where $\phi$ is defined by 
$$\phi(a)\xi(x,y)=a(y)\xi(x,y)$$ 
for $a\in C(\Gamma),$ $\xi\in C(\Omega_{\alpha,\beta}(\Gamma))$ and $(x,y)\in\Omega_{\alpha,\beta}(\Gamma).$ 
Furthermore, we shall assume the existence of an orthonormal basis\index{Basis! Orthonormal} (see Defintion \ref{ONB}) $\{u_i\}_{i=0}^{N-1}$ for $\cE_{\alpha,\beta}(\Gamma)$.\index{$\cO_{\alpha,\beta}(\Gamma)$}

In order to construct our exact sequence for $K_*(\cO_{\alpha,\beta}(\Gamma))$, note the short exact sequence
$$\begindc{0}[10]
\obj(0,0)[O]{$0$}
\obj(5,0)[A]{$\ker q$}
\obj(11,0)[B]{$\cT_{\alpha,\beta}(\Gamma)$}
\obj(18,0)[C]{$\cO_{\alpha,\beta}(\Gamma)$}
\obj(24,0)[D]{$0,$}
\mor{A}{B}{$\iota$}
\mor{B}{C}{$q$}
\mor{C}{D}{$$}
\mor{O}{A}{$$}
\enddc$$
where $q:\cT_{\alpha,\beta}(\Gamma)\to\cO_{\alpha,\beta}(\Gamma)$ is the canonical quotient map and $\iota:\ker q\to \cT_{\alpha,\beta}(\Gamma)$ is the inclusion homomorphism, induces the six-term exact sequence of $K$-groups (see \cite{RLL})

\begin{equation} \label{6term1}
\begindc{0}[10]
\obj(0,2)[A]{$K_0(\mbox{ker } q)$}
\obj(10,2)[B]{$K_0(\cT_{\alpha,\beta}(\Gamma))$}
\obj(20,2)[C]{$K_0(\cO_{\alpha,\beta}(\Gamma))$}
\obj(0,-2)[F]{$K_1(\cO_{\alpha,\beta}(\Gamma))$}
\obj(10,-2)[E]{$K_1(\cT_{\alpha,\beta}(\Gamma))$}
\obj(20,-2)[D]{$K_1(\mbox{ker } q)$}
\mor{A}{B}{$\iota_*$}
\mor{B}{C}{$q_*$}
\mor{C}{D}{$\delta_0$}
\mor{D}{E}{$\iota_*$}
\mor{E}{F}{$q_*$}
\mor{F}{A}{$\delta_1$}
\enddc
\end{equation}\\
\indent Let $(T,\tilde\pi)$ denote the universal Toeplitz representation on $\cE_{\alpha,\beta}(\Gamma);$ that is, $\pi=q\circ\tilde\pi$
is the morphism $C(\Gamma)\to\cO_{\alpha,\beta}(\Gamma).$ As shown in \cite[Theorem 4.4]{Pims}, the homomorphism $\tilde\pi:C(\Gamma)\to\cT_{\alpha,\beta}(\Gamma)$ induces an isomorphism of $K_i(C(\Gamma))$ onto $K_i(\cT_{\alpha,\beta}(\Gamma))$. Thus we may replace $K_i(\cT_{\alpha,\beta}(\Gamma))$ with $K_i(C(\Gamma))$ provided we can identify the resulting maps. We intend to show that (\ref{6term1}) induces the six-term exact sequence

\begin{equation} 
\begindc{0}[10]
\obj(0,2)[A]{$K_0(C(\Gamma))$}
\obj(10,2)[B]{$K_0(C(\Gamma))$}
\obj(20,2)[C]{$K_0(\cO_{\alpha,\beta}(\Gamma))$}
\obj(0,-2)[F]{$K_1(\cO_{\alpha,\beta}(\Gamma))$}
\obj(10,-2)[E]{$K_1(C(\Gamma))$}
\obj(20,-2)[D]{$K_1(C(\Gamma))$}
\mor{A}{B}{$1-\Omega_*$}
\mor{B}{C}{$\pi_*$}
\mor{C}{D}{$$}
\mor{D}{E}{$1-\Omega_*$}
\mor{E}{F}{$\pi_*$}
\mor{F}{A}{$$}
\enddc
\end{equation}\\
\noindent for $\pi=q\circ\tilde\pi:C(\Gamma)\to\cO_{\alpha,\beta}(\Gamma)$ and an appropriately chosen homomorphism $\Omega: C(\Gamma)\to M_N(C(\Gamma)).$

\begin{lemma}\label{Omega}\rm Define $\Omega:C(\Gamma)\to M_N(C(\Gamma))$ by $\Omega(a)=(\innprod{u_i}{a\cdot u_j})_{i,j=0}^{N-1}$. Then $\Omega$ is a unital homomorphism and $\Omega(\alpha^\#(a))$ is the diagonal matrix $\beta^\#(a)1_N$ for all $a\in C(\Gamma).$\\
\pf Let $a,b\in C(\Gamma)$. Then the $(i,j)$-entry of $\Omega(a)\Omega(b)$ is
\begin{align*}
(\Omega(a)\Omega(b))_{i,j}&=\sum_{k=0}^{N-1}\innprod{u_i}{a\cdot u_k}\innprod{u_k}{b\cdot u_j}\\
&=\innprod{u_i}{a\cdot(\sum_k u_k\cdot\innprod{u_k}{b\cdot u_j})}\\ 
&=\innprod{u_i}{a\cdot (b\cdot u_j)}\\
&=\Omega(ab)_{i,j}.
\end{align*}
Furthermore, for $a^*$ denoting the map $a^*(x)=\overline{a(x)}$ for $x\in\Gamma,$
$$\Omega(a^*)=(\innprod{u_i}{a^*\cdot u_j})_{i,j} =(\innprod{a\cdot u_i}{u_j})_{i,j}
=(\innprod{u_j}{a\cdot u_i}^*)_{i,j}=\Omega(a)^*$$
and
$$\Omega(1)=(\innprod{u_i}{u_j})_{i,j}=(\delta_i^j)_{i,j}=1_N.$$
Finally, let $x\in \Gamma$. Then  
\begin{align*}
\Omega(\alpha^\#(a))_{i,j}(x)&=\innprod{u_i}{\alpha^\#(a)\cdot u_j}(x)\\
&=\int_{r^{-1}(x)} \overline{u_i(e)}a(\alpha(s(e)))u_j(e)\,d\lambda_x(e)\\
&=\int_{r^{-1}(x)} \overline{u_i(e)}a(\beta(x))u_j(e)\,d\lambda_x(e)\\
&=a(\beta(x))\int_{r^{-1}(x)} \overline{u_i(e)}u_j(e)\,d\lambda_x(e)\\
&=a(\beta(x))\innprod{u_i}{u_j}(x)\\
&=\delta_i^j\beta^\#(a)(x);
\end{align*}
hence, $\Omega(\alpha^\#(a))=\beta^\#(a)1_N.$\\\sq
\end{lemma}

In order to describe $\ker q$, use the notation $\cE^{\otimes k}:=\cE_{\alpha,\beta}(\Gamma)^{\otimes k}$ 
for the $k$-fold internal tensor product of $C^*$-correspondences (\cite{lan})\index{C$\mbox{}^*$-correspondence! Tensor Product} 
$\cE_{\alpha,\beta}(\Gamma)\otimes\cdots\otimes\cE_{\alpha,\beta}(\Gamma),$ 
which is itself a $C^*$-correspondence over $A=C(\Gamma).$ For the universal covariant representation $(T,\tilde\pi):\cE_{\alpha,\beta}(\Gamma)\to\cT_{\alpha,\beta}(\Gamma)$ (that is, $q(T_j)$ is the isometry $S_j$ with $T(u_j)=T_j,$) there is, in fact, a Toeplitz representation $(T^{\otimes k},\tilde\pi)$ of $\cT_{\alpha,\beta}(\Gamma)$ such that $T^{\otimes k}(\xi)=\prod_{i=1}^k T(\xi_i)$ for all elementary tensors $\xi=\xi_1\otimes\cdots\otimes\xi_k$ where $\xi_i\in\cE_{\alpha,\beta}(\Gamma)$ (see \cite[Proposition 1.8]{FR} where the term ``Hilbert bimodule'' is used instead of $C^*$-correspondence.) 
Note $\cE_{\alpha,\beta}(\Gamma)^{\otimes 0}:=C(\Gamma)$ and $T^{\otimes 0}:=\tilde\pi.$ 
By \cite[Lemma 2.4]{FR},\index{$\cT_{\alpha,\beta}(\Gamma)$}
$$\cT_{\alpha,\beta}(\Gamma)=\overline{\mbox{span}}\{T^{\otimes k}(\xi)T^{\otimes k^\prime}(\eta)^*\st k,k^\prime\ge 0,\xi\in\cE^{\otimes k}, \eta\in\cE^{\otimes k^\prime}\}.$$
Next,  let $p=\sum_{j=0}^{N-1} T_j T_j^*.$ The proceeding lemmas and propositions are essentially those found in \cite[Lemma 3.2, Lemma 3.3 \& Proposition 3.4]{EaHR} with some changes.

\begin{lemma}\rm \label{Lem1} With the preceding notation:
\begin{enumerate} 
\item $p$ is a projection which commutes with $\tilde\pi(a)$ for all $a\in C(\Gamma)$
\item $1-p$ is a full projection in $\ker q$
\item $(1-p)T^{\otimes k}(\xi)=0$ for all $\xi\in\cE^{\otimes k}$ and $k\ge 1$
\item $\ker q=\overline{\mbox{span}}\{T^{\otimes k}(\xi)(1-p)T^{\otimes k^\prime}(\eta)^*\st k,k^\prime\ge0,\xi\in\cE^{\otimes k}, \eta\in\cE^{\otimes k^\prime}\}$
\end{enumerate}
\pf (1) Recall that $T_i^*T_j=\tilde{\pi}(\innprod{u_i}{u_j})=\delta_i^j.$ Thus, $p^2=p$ and $p^*=p.$ Furthermore,
\begin{align*}
p\tilde\pi(a)p&=\sum_{j,k=1}^{N-1}T_jT_j^*\tilde\pi(a)T_kT_k^*=\sum_{j,k=0}^{N-1} T_j\tilde\pi(\innprod{u_j}{a\cdot u_k})T_k^*\\
&=\sum_{j,k=0}^{N-1} T(u_j\innprod{u_j}{a\cdot u_k})T_k^*=\sum_{k=0}^{N-1}T(\sum_{j=0}^{N-1}u_j\innprod{u_j}{a\cdot u_k})T_k^*\\
&=\sum_{k=0}^{N-1}T(a\cdot u_k)T_k^*=\tilde\pi(a)p
\end{align*}
and so,
$$p\tilde\pi(a)=(\tilde\pi(a)^*p)^*=(p\tilde\pi(a)^*p)^*=p\tilde\pi(a)p=\tilde\pi(a)p.$$

(2) Recall $\phi(a)=\sum_{j=0}^{N-1}\theta_{a\cdot u_j,u_j}$, so
$$\Psi_T(\phi(a))=\sum_{j=0}^{N-1}T(a\cdot u_j)T(u_j)^*=\tilde\pi(a)p$$
and $$q(1-p)=q(\tilde\pi(1)-\tilde\pi(1)p)=q(\tilde\pi(1)-\Psi_T(\phi(1)))=0.$$ 
Hence, $1-p=\tilde\pi(1)-\tilde\pi(1)p\in\ker q$ and since $\ker q$ is the ideal in $\cT_{\alpha,\beta}(\Gamma)$ generated by $\{\tilde\pi(a)-\Psi_T(\phi(a))\st a\in C(\Gamma)\}$ and $1-p\in\ker q$, $\ker q$ is the ideal generated by $\{\tilde\pi(a)(1-p)\st a\in C(\Gamma)\}.$ Hence $1-p$ is full.

(3) Let $\xi\in\cE_{\alpha,\beta}(\Gamma)$ then 
$$pT(\xi)=\sum_{j=0}^{N-1} T(u_j)T(u_j)^*T(\xi)=\sum_{j=0}^{N-1}T(u_j\innprod{u_j}{\xi})=T(\xi).$$
Thus, $(1-p)T(\xi)=0$. Now for $k>1$, let $\xi=\xi_1\otimes...\otimes\xi_k$. Then
$$(1-p)T^{\otimes k}(\xi)=(1-p)\prod_{j=0}^kT(\xi_j)=0$$
and hence, by linearity and continuity, (3) has been proven.

(4) Since $\ker q=\cT_{\alpha,\beta}(\Gamma)(1-p)\cT_{\alpha,\beta}(\Gamma),$ the description of $\cT_{\alpha,\beta}(\Gamma)$ preceding Lemma \ref{Lem1} paired with (3) gives the desired result.\\\sq
\end{lemma}

\begin{lemma}\label{Rho}\rm There exists a homomorphism $\rho:C(\Gamma)\to\ker q\subset\cT_{\alpha,\beta}(\Gamma)$ 
such that $\rho(a)=\tilde\pi(a)(1-p)$ and $\rho$ is an isomorphism of $C(\Gamma)$ onto the full corner $C^*$-algebra 
$(1-p)\mbox{ker }q(1-p)$.\\
\pf By the previous lemma,
$$(1-p)\tilde\pi(a)(1-p)=\tilde\pi(a)(1-p)\in\ker q.$$
Thus, $\rho(a)=\tilde\pi(a)(1-p)$ defines a homomorphism $\rho:C(\Gamma)\to (1-p)\ker q(1-p)\subset\ker q.$
Using the previous lemma,
\begin{align*}
(1-p)\ker q(1-p)&=\overline{\mbox{span}}\{(1-p)T^{\otimes k}(\xi)(1-p)T^{\otimes k^\prime}(\eta)^*(1-p)\st k,k^\prime\ge0, \xi\in\cE^{\otimes k},\eta\in\cE^{\otimes k^\prime}\}\\
&=\overline{\mbox{span}}\{(1-p)\tilde\pi(a)(1-p)\tilde\pi(b)^*(1-p)\st a,b\in C(\Gamma)\}\\
&=\overline{\mbox{span}}\{\tilde\pi(a)(1-p)\st a\in C(\Gamma)\}=\mbox{ran }\rho.
\end{align*}
Hence, $\rho$ is surjective. 

In order to show the injectivity of $\rho$, choose a faithful representation $\pi_0:C(\Gamma)\to B(H)$ and consider
the Fock representation $(T_F,\pi_F)$ of $\cE_{\alpha,\beta}(\Gamma)$ \index{Fock Representation}
induced from $\pi_0$ as described in \cite[Example 1.4]{FR}. The underlying space of this Fock representation is 
$F(\cE_{\alpha,\beta}(\Gamma))\otimes_A H:=\oplus_{k\ge0}(\cE^{\otimes k}\otimes_A H)$
where $A=C(\Gamma)$ acts diagonally on the left and $\cE_{\alpha,\beta}(\Gamma)$ acts by creation operators.
\index{Creation Operator} 
Then $T_F(\xi)^*$ is an annihilation operator\index{Annihilation Operator} vanishing on the subspace $A\otimes_A H$ of 
$F(\cE_{\alpha,\beta}(\Gamma))\otimes_A H.$
Now, for $a\in A,$ 
$$0=(T_F\times\pi_F)(\rho(a))=(T_F\times\pi_F)(\tilde\pi(a)(1-p)=\pi_F(a)(1-\sum_{j=0}^{N-1}T_F(u_j)T_F(u_j)^*).$$
Since $T_F(u_j)^*$ vanishes on $A\otimes_A H$, we have that $\rho(a)=0$ implies 
$$\pi_F(a)(1-\sum_{j=0}^{N-1}T_F(u_j)T_F(u_j)^*)(1\otimes_A h)=0$$ 
for all $h\in H$ and so, $\pi_F(a)(1\otimes_A h)=0$ for all $h\in H.$ Thus, $a\otimes_A h=0$ for all $h\in H$ and hence, $\pi_0(a)h=0$ for all $h\in H$ which implies $a=0$ since $\pi_0$ is faithful. Hence, $\rho$ is injective.\\\sq
\end{lemma}

\begin{lemma}\cite[Lemma 3.5]{EaHR} \label{KLem1} \rm 
Suppose that $A$ is a $C^*$-algebra, $r\ge1$ and $N\ge 2$ are integers, and
$$\{b_{j,s;k,t}\st  0\le j,k < N\mbox{ and } 0\le s,t < r\}$$
is a subset of $A.$ For $m,n$ satisfying $0\le m,n<rN-1,$ define
$$c_{m,n}=b_{j,s;k,t}\mbox{ where }m=sN+j\mbox{ and }n=lN+k,\mbox{ and}$$
$$d_{m,n}=b_{j,s;k,t}\mbox{ where }m=jr+s\mbox{ and }n=kr+t.$$
Then there is a scalar unitary permutation matrix $U$ such that the matrices $C:=(c_{m,n})_{m,n}$ and $D_{m,n}:=(d_{m,n})_{m,n}$ are related by $C=UDU^*.$ 
\end{lemma}

The following is standard (and also appears in \cite{EaHR}):

\begin{lemma}\label{KLem2} \rm 
Suppose that $S$ is an isometry in a unital $C^*$-algebra $A$. Then 
$$U:=\mat{S}{1-SS^*}{0}{S^*}$$
is a unitary element of $M_2(A)$ and its class in $K_1(A)$ is the identity.
\end{lemma}

\begin{proposition}\rm  Let $(T,\tilde\pi)$ denote the universal Toeplitz representation on $\cE_{\alpha,\beta}(\Gamma)$ and
let $\{u_j\}_{j=0}^{N-1}$ be an  orthonormal basis of $\cE_{\alpha,\beta}(\Gamma).$
Further, let $p=\sum_{j=0}^{N-1} T_j T_j^*$ where $T_j=T(u_j).$ Then, 
with the maps $\Omega:C(\Gamma)\to M_N(C(\Gamma))$ and $\rho:C(\Gamma)\to\ker q\subset\cT_{\alpha,\beta}(\Gamma)$
defined by 
$$\Omega(a)=(\innprod{u_i}{a\cdot u_j})_{i,j=0}^{N-1}$$
and 
$$\rho(a)=\tilde\pi(a)(1-p)$$
as in Lemmas \ref{Omega} and \ref{Rho},
the following two diagrams $(i=0,1)$ commute:
\begin{equation}
\begindc{0}[10]
\obj(0,2)[A]{$K_i(C(\Gamma))$}
\obj(10,2)[B]{$K_i(C(\Gamma))$}
\obj(0,-2)[D]{$K_i(\mbox{ker }q)$}
\obj(10,-2)[C]{$K_i(\cT_{\alpha,\beta}(\Gamma))$}
\mor{A}{B}{$1-\Omega_*$}
\mor{B}{C}{$\tilde\pi_*$}
\mor{D}{C}{$\iota_*$}
\mor{A}{D}{$\rho_*$}
\enddc
\end{equation}
\pf First, let $i=0.$ Let $z=(z_{s,t})\in M_r(C(\Gamma))$ be a projection and let $\tilde\pi_r$ denote the augmentation map,
$\tilde\pi\otimes\mbox{id}_r,$ of $\tilde\pi$ on $M_r(C(\Gamma)).$
Then 
$$\rho_*([z])=[(\rho(z_{s,t}))_{s,t}]=[(\tilde\pi(z_{s,t})(1-p))_{s,t}]=[\tilde\pi_r(z)]-[\tilde\pi_r(z)(p1_r)]$$
and
$$\tilde\pi_*\circ(1-\Omega_*)([z])=[\tilde\pi_r(z)]-\tilde\pi_*\circ\Omega_*([z]).$$
Hence, it suffices to show that
$$\tilde\pi_*\circ\Omega_*([z])=[\tilde\pi_r(z)(p1_r)].$$
Note that
$$\Omega_*([z])=[(\Omega(z_{s,t}))_{s,t}]=[((\innprod{u_j}{z_{s,t}\cdot u_k})_{j,k})_{s,t}],$$
so
$$\tilde\pi_*\circ\Omega_*([z])=[((\tilde\pi(\innprod{u_j}{z_{s,t}\cdot u_k}))_{j,k})_{s,t}]=[\tilde\pi_{rN}\circ\Omega_r(z)].$$
Set $b_{j,s;k,t}=\tilde\pi(\innprod{u_j}{z_{s,t}\cdot u_k})$ and $C=(c_{m,n})_{m,n}=\tilde\pi_{rN}(\Omega_r(z))$ as in
Lemma \ref{KLem1}.

Let $$T=\begin{pmatrix}
T_01_r &T_11_r&...&T_{N-1}1_r\\
0&0&...&0\\
\vdots&\vdots&...&\vdots\\
0&0&...&0
\end{pmatrix}\in M_N(M_r(\cO_{\alpha,\beta}(\Gamma))).$$
Then $TT^*=p1_r\oplus 0_{r(N_1)}$ and since $\tilde\pi_r(z)$ is a projection which commutes with $p1_r,$
$$(\tilde\pi_r(z)\oplus 0_{r(N-1)})T$$
is a partial isometry which implements a Murray-von Neumann equivalence between
$$T^*(\tilde\pi_r(z)\oplus 0_{r(N-1)})T$$
and
$$(\tilde\pi_r(z)\oplus 0_{r(N-1)})TT^*(\tilde\pi_r(z)\oplus 0_{r(N-1)})=\tilde\pi_r(z)(p1_r)\oplus 0_{r(N-1)};$$
thus,
$$[\tilde\pi_r(z)(p1_r)]=[T^*(\tilde\pi_r(z)\oplus 0_{r(N-1)})T].$$
Furthermore,
$$T^*(\tilde\pi_r(z)\oplus 0_{r(N-1)})T=T^*\begin{pmatrix} 
\tilde\pi_r(z)T_0&...&\tilde\pi_r(z)T_{N-1}\\
0&...&0\\
\vdots&...&\vdots
\end{pmatrix} =(T_j^*\tilde\pi_r(z)T_k)_{j,k}$$
so the $(j,k)$ entry is $(\tilde\pi(\innprod{u_j}{z_{s,t}\cdot u_k}))_{s,t}.$
Recall $b_{j,s;k,t}=\tilde\pi(\innprod{u_j}{z_{s,t}\cdot u_k})$ and so
$T^*(\tilde\pi_r(z)\oplus 0_{r(N-1)})T=D=(d_{m,n})_{m,n}$ as in Lemma \ref{KLem1}. Thus, by Lemma \ref{KLem1}, there
exists a unitary $U$ such that $C=U^*DU$ which gives us
$$[\tilde\pi_r(z)(p1_r)]=[D]=[C]=[\tilde\pi_{rN}\circ\Omega_R(z)]$$ 
as desired.

For the case $i=1,$ let $u\in M_r(C(\Gamma))$ be a unitary. Note $\rho_*: K_1(C(\Gamma))\to K_1(\ker q)$ is the composition of a unital isomorphism
of $C(\Gamma)$ onto $(1-p)\ker q(1-p)$ with the inclusion of $(1-p)\ker q(1-p)$ as a full corner in the non-unital algebra $\ker q;$
that is, $[u]\mapsto [\rho_r(u)]=[\tilde\pi_r(u)((1-p)1_r)]\mapsto [\tilde\pi_r(u)((1-p)1_r)+p1_r]\in K_1((\ker q)^+)=K_1(\ker q).$
Furthermore,
$$\tilde\pi_*\circ\Omega_*([u])=[\tilde\pi_r(u)]-[\tilde\pi_{rN}\circ\Omega_r(u)]$$
and hence, we need only show
$$[(\tilde\pi_r(u)((1-p)1_r)+p1_r)\oplus 1_{r(N-1)}]=[\tilde\pi_r(u)\oplus 1_{r(N-1)}]-[\tilde\pi_{rN}\circ\Omega_r(u)]$$
in $K_1(\cT_{\alpha,\beta}(G)).$

We take a brief moment to make an aside: If $C\in M_{2rN}(\cT_{\alpha,\beta}(\Gamma))$ is invertible with $K_1$-class the identity $1$ then the $K_1$-class is unchanged by pre- and post-multiplication by $C$.
In particular, when $C$ is equal to:
\begin{enumerate}
\item (Lemma \ref{KLem2}) a unitary of the form
$$\mat{S}{1-SS^*}{0}{S^*}$$ where $S$ is an isometry
\item an upper- or lower-triangular matrix of the form
$$\mat{1}{A}{0}{1}\qquad\mbox{or}\qquad\mat{1}{0}{A}{1}$$
(which are connected to $1$ via $t\mapsto\mat{1}{tA}{0}{1}$ and likewise for the transpose)
\item any constant invertible matrix in $\mbox{GL}_{2rN}(\zC)$ (because $\mbox{GL}_{2rN}(\zC)$ is path connected);
this implies that row and column operations may be used without changing the $K_1$-class. 
\end{enumerate} 
$$\mbox{Recall: }T=\begin{pmatrix}
T_01_r &T_11_r&...&T_{N-1}1_r\\
0&0&...&0\\
\vdots&\vdots&...&\vdots\\
0&0&...&0
\end{pmatrix}.$$
With this in mind, calculate
\begin{align*}
[(\tilde\pi_r(u)&((1-p)1_r)+p1_r)\oplus 1_{r(N-1)}]\\
&=\Big[\mat{(\tilde\pi_r(u)((1-p)1_r)+p1_r)\oplus 1_{r(N-1)}}{0_{rN}}{0_{rN}}{1_{rN}}\Big]
\Big[\mat{T}{1_{rN}-TT^*}{0_{rN}}{T^*}\Big]\\
&=\Big[\mat{(\tilde\pi_r(u)((1-p)1_r)+p1_r)\oplus 1_{r(N-1)}}{0_{rN}}{0_{rN}}{1_{rN}}\Big]
\Big[\mat{T}{(1-p)1_r\oplus 1_{r(N-1)}}{0_{rN}}{T^*}\Big]\\
&=\Big[\mat{((\tilde\pi_r(u)(1-p)1_r)+p1_r)\oplus 1_{r(N-1)})T}{\tilde\pi(u)((1-p)1_r)\oplus 1_{r(N-1)}}{0_{rN}}{T^*}\Big]
\end{align*}
and recall $(1-p)T_i=0$ by Lemma \ref{Lem1}(3), hence $(1-p)1_rT=0$ and
\begin{align*}
[(\tilde\pi_r(u)&((1-p)1_r)+p1_r)\oplus 1_{r(N-1)}]\\
&=\Big[\mat{T}{\tilde\pi_r(u)((1-p)1_r)\oplus 1_{r(N-1)}}{0_{rN}}{T^*}\Big]\\
&=\Big[\mat{T}{\tilde\pi_r(u)((1-p)1_r)\oplus 1_{r(N-1)}}{0_{rN}}{T^*}\Big]
\Big[\mat{1_{rN}}{T^*(\tilde\pi_r(u)\oplus 1_{r(N-1)})}{0_{rN}}{1_{rN}}\Big]\\
&=\Big[\mat{T}{\tilde\pi_r(u)\oplus 1_{r(N-1)}}{0_{rN}}{T^*}\Big]
\end{align*}
since $TT^*=p1_r\oplus 0_{r(N-1)}$ and $(p1_r)\tilde\pi_r(u)=\tilde\pi_r(u)(p1_r).$ Using elementary operations, compute
\begin{align*}
[(\tilde\pi_r(u)&((1-p)1_r)+p1_r)\oplus 1_{r(N-1)}]\\
&=\Big[\mat{\tilde\pi_r(u)\oplus 1_{r(N-1)}}{T}{T^*}{0_{rN}}\Big]\\
&=\Big[\mat{\tilde\pi_r(u)\oplus 1_{r(N-1)}}{T}{T^*}{0_{rN}}\Big]\Big[\mat{1_{rN}}{-(\tilde\pi_r(u*)\oplus 1_{r(N-1)})T}{0_{rN}}{1_{rN}}\Big]\\
&=\Big[\mat{\tilde\pi_r(u)\oplus 1_{r(N-1)}}{0_{rN}}{T^*}{-T^*(\tilde\pi_r(u^*)\oplus 1_{r(N-1)})T}\Big]\\
&=\Big[\mat{1_{rN}}{0_{rN}}{-T^*(\tilde\pi_r(u^*)\oplus 1_{r(N-1)})}{1_{rN}}\Big]\Big[\mat{\tilde\pi_r(u)\oplus 1_{r(N-1)}}{0_{rN}}{T^*}{-T^*(\tilde\pi_r(u^*)\oplus 1_{r(N-1)})T}\Big]\\
&=\Big[\mat{\tilde\pi_r(u)\oplus 1_{r(N-1)}}{0_{rN}}{0_{rN}}{-T^*(\tilde\pi_r(u^*)\oplus 1_{r(N-1)})T}\Big]\Big[\mat{1_{rN}}{0_{rN}}{0_{rN}}{-1_{rN}}\Big]\\
&=[\tilde\pi_r(u)\oplus 1_{r(N-1)}]+[T^*(\tilde\pi_r(u^*)\oplus 1_{r(N-1)})T].
\end{align*}
Furthermore,
$$[T^*(\tilde\pi_r(u^*)\oplus 1_{r(N-1)})T]=[\tilde\pi_{rN}(\Omega_r(u^{-1}))]=-[\tilde\pi_{rN}(\Omega_r(u))].$$
Hence,
$$[(\tilde\pi_r(u)((1-p)1_r)+p1_r)\oplus 1_{r(N-1)}]=[\tilde\pi_r(u)\oplus 1_{r(N-1)}]-[\tilde\pi_{rN}\circ\Omega_r(u)]$$
as desired.\\ \sq
\end{proposition}

\begin{theorem}\label{NewSixTerm}\rm  \index{$\cO_{\alpha,\beta}(\Gamma)$! Six Term Exact Sequence}
Let $(S,\pi)=q\circ(T,\tilde\pi)$ be the universal Cuntz-Pimsner covariant representation of $\cE_{\alpha,\beta}(\Gamma)$ in 
$\cO_{\alpha,\beta}(\Gamma).$ Then the following diagram is exact:
\begin{equation}\label{6term2} 
\begindc{0}[10]
\obj(0,2)[A]{$K_0(C(\Gamma))$}
\obj(10,2)[B]{$K_0(C(\Gamma))$}
\obj(20,2)[C]{$K_0(\cO_{\alpha,\beta}(\Gamma))$}
\obj(0,-2)[F]{$K_1(\cO_{\alpha,\beta}(\Gamma))$}
\obj(10,-2)[E]{$K_1(C(\Gamma))$}
\obj(20,-2)[D]{$K_1(C(\Gamma))$}
\mor{A}{B}{$1-\Omega_*$}
\mor{B}{C}{$\pi_*$}
\mor{C}{D}{$\rho_*^{-1}\circ\delta_1$}
\mor{D}{E}{$1-\Omega_*$}
\mor{E}{F}{$\pi_*$}
\mor{F}{A}{$\rho_*^{-1}\circ\delta_0$}
\enddc
\end{equation}
\pf  Note $\rho:C(\Gamma)\to\mbox{ker }q$ is an isomorphism onto a full corner, implying $\rho_*$ is an isomorphism. Further note
$\tilde\pi_*:K_i(C(\Gamma))\to K_i(\cT_{\alpha,\beta}(\Gamma))$ is an isomorphism (see comments prior to Lemma \ref{Omega}). Then (\ref{6term1}) and the previous proposition give the stated result.\\\sq
\end{theorem}

\section{$K$-groups of $\cO_{F,G}(\Gamma)$}In this section, the approach of \cite{EaHR} is made easier and extended.
For this section, let $\alpha,\beta\in\mbox{End}(C(\zT^d))$ defined by $\alpha=\sigma_F^\#$ and $\beta=\sigma_G^\#$ where $F=\diag{a_1,...,a_d}\in M_d(\zZ)^+$ and $G\in M_d(\zZ)$ such that $\det F>0$ and $\det G\ne0$. We know there exists an orthonormal basis for $\cE_{F,G}(\zT^d),$ $\{u_j\st 0\le j\le N-1\};$ this is the basis $\{u_\nu\}_{\nu\in\fI(F)},$ described in Section 3.3, reindexed by $0,1,...,N-1.$ Let $U_j$ be the unitary defined by $U_j(x)=x_j$ for $x=(x_i)_{i=1}^d\in\zT^d$ for $j\in\{1,...,d\}.$ Further,
let 
$$\fI_k=\{J\subset\{1,...,d\}\st \abs{J}=k, J=\{j_1<...<j_k\}\}$$ 
and $J^\prime=\{1,...,d\}\setminus J$ in increasing order. Define
$$\mathfrak E_k=\begin{cases} 
\{[1]_0\}&\mbox{if }k=0\\ 
\{[U_J]_0=[U_{j_1}]_0\wedge...\wedge[U_{j_k}]_0\st J\in\fI_k\}&\mbox{if $k>0$ is even}\\
\{[U_J]_1=[U_{j_1}]_1\wedge...\wedge[U_{j_k}]_1\st J\in\fI_k\}&\mbox{if $k>0$ is odd}
\end{cases}$$ If it is understood, the notation $[\cdot]$ will be used in lieu of $[\cdot]_i.$  

It is well known (see \cite{Ji} and \cite[Example 3.11 and 3.15]{Hat}) that 
$$K_0(C(\zT^d))\cong\bigwedge_{\rm evens}\zZ^d=\zZ^{2^{d-1}}$$ 
with basis $\{\mathfrak E_k\}_{k \mbox{ even}}$ and 
$$K_1(C(\zT^d))\cong\bigwedge_{\rm odds}\zZ^d=\zZ^{2^{d-1}}$$
with basis $\{\mathfrak E_k\}_{k \mbox{ odd}}.$ 
For subsets $J$ and $I$ of the same size, define $F_{J,I}$ to be the square submatrix of $F$ whose 
entries belong to the rows in $J$ and the columns in $I$. 

With these identifications, the ($K_1$-group) induced map $\alpha_*\vert_{\bigwedge^1\zZ^d}:\mbox{span}\{[U_j]\}_{j=1}^d\to\mbox{span}\{[U_j]\}_{j=1}^d$
is multiplication by $F^T=F,$ and $\beta_*\vert_{\bigwedge^1\zZ^d}:\mbox{span}\{[U_j]\}_{j=1}^d\to\mbox{span}\{[U_j]\}_{j=1}^d$
is multiplication by $G^T,$ the transpose of $G.$ We have 
\begin{align*}
\beta_*([U_j])&=[\beta(U_j)]=[U^{G_j}]=[\prod_k U_k^{b_{jk}}]=\sum_k b_{jk}[U_k].
\end{align*}
Do likewise to prove $\alpha_*\vert_{\bigwedge^1\zZ^d}$ is multiplication by $F.$
One can also check that $\alpha_*$ and $\beta_*$ act on $\bigwedge^0\zZ^d$ by 
$$\alpha_*[1]=[\alpha(1)]=1=[\beta(1)]=\beta_*[1]$$ 
since $\alpha$ and $\beta$ are group homomorphisms. 

\begin{lemma}\label{Bs}\rm For $1\le k\le d$, the matrix $A_k$ representating $\alpha_*\vert :\bigwedge^k\zZ^d\to\bigwedge^k\zZ^d$ with respect to the basis $\mathfrak E_k$ is the diagonal matrix $\diag{a_I}_{I\in\fI_k}$ ($a_I=\prod_{i\in I}a_i$) and matrix $B_k$ representating $\beta_*\vert:\bigwedge^k\zZ^d\to\bigwedge^k\zZ^d$ is $(\det G_{JI})_{I,J\in\fI_k}$.\\
\pf Begin by noting, $A_k=\bigwedge^kF^T$ and $B_k=\bigwedge^kG^T.$
Let $[U_I]\in\mathfrak E_k$ with $I=\{i_1<...<i_k\}$. Then
\begin{align*} 
\beta_*([U_I])&=(\bigwedge^k G^T)[U_I]=(\bigwedge^k G^T)([U_{i_1}]\wedge...\wedge[U_{i_k}])\\
&=G^T[U_{i_1}]\wedge ...\wedge G^T[U_{i_k}]\\
&=\sum_{m_1,...,m_k=1}^d b_{i_1,m_1}...b_{i_k,m_k}([U_{m_1}]\wedge ...\wedge [U_{m_k}])\\
&=\sum_{[U_J]\in\mathfrak E_k, J=\{m_1,...,m_k\}}  b_{i_1,m_1}...b_{i_k,m_k}[U_J]\\
&=\sum_{[U_J]\in\mathfrak E_k}\sum_{\sigma\in S_k}  b_{i_1,\sigma(j_1)}...b_{i_k,\sigma(j_k)}([U_{\sigma(j_1)}]\wedge ... \wedge [U_{\sigma(j_k)}])\\
&=\sum_{[U_J]\in\mathfrak E_k}\sum_{\sigma\in S_k}  (-1)^{\deg \sigma}b_{i_1,\sigma(j_1)}...b_{i_k,\sigma(j_k)}([U_{j_1}]\wedge ... \wedge [U_{j_k}])\\
&=\sum_{[U_J]\in\mathfrak E_k} \det G_{I,J}[U_J]
\end{align*}
where $S_k$ denotes the symmetric group on $k$ elements. The result for $A_k$ follows by specializing to the diagonal case.\\\sq
\end{lemma}

Let $A_k$ and $B_k$ (for $k=0,...,d$) denote the matrices described in Lemma \ref{Bs}; that is, $A_0=B_0=1,$ $A_k=\diag{a_I}_{I\in\fI_k}$ and $B_k=(\det G_{JI})_{I,J\in\fI_k}$ for $k\in\{1,...,d\}.$ From Lemma \ref{Omega}, the map $\Omega(\sigma_F^\#):C(\zT^d)\to M_N(C(\zT^d))$ is the diagonal matrix  $\Omega(\sigma_F^\#(a))=\sigma_G^\#(a)1_N$ and so, for $\alpha=\sigma_F^\#,$ $\beta=\sigma_G^\#$ and $d=(d_{s,t})\in M_r(C(\zT^d)),$
\begin{align*}
\Omega_*\circ\alpha_*([d])&=[(\Omega\circ\alpha)_r(d)]\\
&=[(\Omega(\alpha(d_{s,t}))_{s,t}]\\
&=[(\beta(d_{s,t})1_N)_{s,t}]\\
&=N[\beta_r(d)]\\
&=N\beta_*[d]
\end{align*}
where $\Omega_r,\alpha_r,\beta_r$ denote the appropriate augmented maps on $M_r(C(\zT^d))$. Using the previous lemma and the above equation, the matrix $C_k$ representing $\Omega_*$ on $\bigwedge^k\zZ^d$ satisfies 
$$C_kA_k=NB_k$$
where it is expected that $C_k$ is a matrix with integer entries for each $k=0,...,d.$

Recall  $A_k=\diag{a_I}_{I\in\fI_k}.$
Let $I^\prime=\{1,...,d\}\setminus I$ ordered so that
$I=\{i_1<...<i_k\}, I^\prime=\{i_{k+1}<...<i_d\}.$  
Then set $C_0=N=\det F\in\zZ$ and
$$C_k=B_k\diag{a_{I^\prime}}_{I\in\fI_k}=(a_{J^\prime}\det G_{JI})_{I,J\in\fI_k}$$
for $k=1,...,d.$ Note that $A_0=B_0=1$ by the calculations before Lemma \ref{Bs}. Hence, 
$$C_0A_0=NB_0$$
and
\begin{align*}
C_kA_k&=B_k\diag{a_{I^\prime}}_{I\in\fI_k}\diag{a_I}_{I\in\fI_k}\\
&=B_k\diag{a_{I\cup I^\prime}}_{I\in\fI_k}\\
&=B_k\diag{a_{1,...,d}}_{I\in\fI_k}\\
&=B_k\diag{N}_{I\in\fI_k}=NB_k
\end{align*} for $k=1,...,d.$

In order to calculate the $K$-theory of $\cO_{F,G}(\zT^d)$, one needs only to calculate 
$\ker(1-C_k)$ and $\mbox{coker}(1-C_k)$ for each $k=0,...,d$ as the next theorem will demonstrate.

\begin{theorem}\label{Alg}\rm \index{$\cO_{F,G}(\zT^d)$! $K$-groups}
Let $F=\diag{a_1,...,a_d}, G\in M_d(\zZ)$ such that $\det G\ne0$ and $a_j\in\zN$ for each $j=1,...,d.$ 
With $C_k$ defined as above and $\cO_{F,G}(\zT^d)$
defined as in Example \ref{ToriQuiver},
\begin{enumerate}
\item $K_0(\cO_{F,G}(\zT^d))=\ds\big(\bigoplus_{0\le k\le d,\,\rm even} \mbox{coker}(1-C_k)\big)
\oplus\big(\bigoplus_{0\le k\le d, \,\rm odd} \ker(1-C_k)\big)$, and
\item$K_1(\cO_{F,G}(\zT^d))=\ds\big(\bigoplus_{0\le k\le d, \,\rm odd} \mbox{coker}(1-C_k)\big)
\oplus\big(\bigoplus_{0\le k\le d, \,\rm even} \ker(1-C_k)\big).$
\end{enumerate}
\pf By \ref{6term2} in Theorem \ref{NewSixTerm},
$$\begindc{0}[15]
\obj(0,2)[A]{$K_0(C(\zT^d)$}
\obj(10,2)[B]{$K_0(C(\zT^d))$}
\obj(20,2)[C]{$K_0(\cO_{F,G}(\zT^d))$}
\obj(0,-2)[F]{$K_1(\cO_{F,G}(\zT^d))$}
\obj(10,-2)[E]{$K_1(C(\zT^d))$}
\obj(20,-2)[D]{$K_1(C(\zT^d))$}
\mor{A}{B}{$\ds\bigoplus_{0\le k\le d, \,\rm even}1-C_k$}
\mor{B}{C}{$$}
\mor{C}{D}{$$}
\mor{D}{E}{$\ds\bigoplus_{0\le k\le d, \,\rm odd}1-C_k$}
\mor{E}{F}{$$}
\mor{F}{A}{$$}
\enddc$$
is exact, so there are two exact sequences 
$$\begindc{0}[20]
\obj(0,0)[O]{$0$}
\obj(5,0)[A]{$\ds\bigoplus_{0\le k\le d, \,\rm even}\mbox{coker}(1-C_k)$}
\obj(11,0)[B]{$K_0(\cO_{F,G}(\zT^d))$}
\obj(17,0)[C]{$\ds\bigoplus_{0\le k\le d, \,\rm odd}\ker(1-C_k)$}
\obj(21,0)[D]{$0$}
\mor{A}{B}{$$}
\mor{B}{C}{$$}
\mor{C}{D}{$$}
\mor{O}{A}{$$}
\enddc$$
and
$$\begindc{0}[20]
\obj(0,0)[O]{$0$}
\obj(5,0)[A]{$\ds\bigoplus_{0\le k\le d, \,\rm odd}\mbox{coker}(1-C_k)$}
\obj(11,0)[B]{$K_1(\cO_{F,G}(\zT^d))$}
\obj(17,0)[C]{$\ds\bigoplus_{0\le k\le d, \,\rm even}\ker(1-C_k)$}
\obj(21,0)[D]{$0$}
\mor{A}{B}{$$}
\mor{B}{C}{$$}
\mor{C}{D}{$$}
\mor{O}{A}{$$}
\enddc$$
which split since $\bigwedge^k\zZ^d$ and hence, $\ker(1-C_k)$ is free for each $k$. Thus, 
$$K_0(\cO_{F,G}(\zT^d))=\big(\bigoplus_{0\le k\le d, \,\rm even} \mbox{coker}(1-C_k)\big)
\oplus\big(\bigoplus_{0\le k\le d, \,\rm odd} \ker(1-C_k)\big)$$ and
$$K_1(\cO_{F,G}(\zT^d))=\big(\bigoplus_{0\le k\le d, \,\rm odd} \mbox{coker}(1-C_k)\big)
\oplus\big(\bigoplus_{0\le k\le d, \,\rm even} \ker(1-C_k)\big).$$
\end{theorem} 

\begin{definition}\rm A matrix $Z\in M_d(\zZ)$ is called an \emph{integer dilation matrix}\index{Integer Dilation Matrix} 
provided each eigenvalue $\lambda$ of $Z$ satisfies $\abs{\lambda}>1.$
\end{definition}

\begin{remark}\rm The case where $F$ is an integer dilation and $G=1_d$ was computed in \cite[Theorem 4.9]{EaHR} where it
was found that 
$$K_0(\cO_{F,1}(\zT^d))=\big(\bigoplus_{0\le k\le d, \,\rm even} \mbox{coker}(1-Q_k)\big)
\oplus\big(\bigoplus_{0\le k\le d, \,\rm odd} \ker(1-Q_k)\big)$$ and
$$K_1(\cO_{F,1}(\zT^d))=\big(\bigoplus_{0\le k\le d, \,\rm odd} \mbox{coker}(1-Q_k)\big)
\oplus\big(\bigoplus_{0\le k\le d, \,\rm even} \ker(1-Q_k)\big)$$
for the matrix $Q_k$ satisfying the relation
$$Q_k(\det F_{JI})_{I,J\in\fI_k}=N1_{d\choose k}.$$ 
Recall the Smith normal form of $F,$ $F=UDV$ where $U,V\in M_d(\zZ)$ are unimodular matrices and 
$D=\diag{a_j}_{j=1}^d\in M_d(\zZ)$ is a positive diagonal matrix. Then by properties of matrix minors,
$$(\det F_{JI})_{I,J\in\fI_k}=U_kD_kV_k$$
where $U_k=(\det U_{JI})_{I,J\in\fI_k},$ $V_k=(\det V_{JI})_{I,J\in\fI_k}$ and $D_k=\diag{a_I}_{I\in\fI_k}.$
Also note for $(U^{-1})_k=(\det (U^{-1})_{JI})_{I,J\in\fI_k},$
$$U_k(U^{-1})_k=(UU^{-1})_k=1_{d\choose k};$$
that is, $U_k$ is unimodular.
Hence, for $G=U^{-1}V^{-1},$ $Q_kU_kD_kV_k=N1_{d\choose k}$ implies
$$U_k^{-1}Q_kU_k=U_k^{-1}V_k^{-1}\diag{a_{I^\prime}}_{I\in\fI_k}=B_k\diag{a_{I^\prime}}_{I\in\fI_k}=C_k.$$
Therefore,
$$\ker(1-C_k)=\ker(U_k^{-1}(1-Q_k)U_k)\cong\ker(1-Q_k)$$
and likewise,
$$\mbox{coker}(1-C_k)\cong\mbox{coker}(1-Q_k);$$
that is, Theorem \ref{Alg} extends Theorem 4.9 of \cite{EaHR}.
\end{remark} 

Here we first consider the case where $F=n1_d$ and $G\in M_d(\zZ).$  To this end, note that $C_k$ has a simple form. First of all, note $A_k$ is $n^k1_{d\choose k}$ and hence,
$$C_k=n^{d-k}B_k. \qquad (\mbox{see the calculations preceding Theorem \ref{Alg}.})$$

\begin{theorem}\label{Alg2} \rm For $d\in\zN$, let $G\in M_d(\zZ)$ $(\det G\ne0)$ and $n\in\zN$. Then
with $C_k=n^{d-k}B_k$ where $B_k=(\det G_{J,I})_{I,J\in\fI_k}\in M_{d\choose k}(\zZ)$ as above, 
\begin{enumerate}
\item If $d>1$ and $n>1$, then 
\begin{enumerate}
\item $K_0(\cO_{n,G}(\zT^d))=
\begin{cases}
\ds\bigoplus_{0\le k\le d,\;\mbox{even}}\mbox{coker}(1-C_k)&\mbox{if $\det G\ne1$}\\
\zZ\oplus\big(\ds\bigoplus_{0\le k\le d-1,\;\mbox{even}}\mbox{coker}(1-C_k)\big)&\mbox{if $\det G=1$}
\end{cases}$
\item $K_1(\cO_{n,G}(\zT^d))=
\begin{cases}
\ds\bigoplus_{0\le k\le d,\;\mbox{odd}}\mbox{coker}(1-C_k)&\mbox{if $\det G\ne1$}\\
\zZ\oplus\big(\ds\bigoplus_{0\le k\le d-1,\;\mbox{odd}}\mbox{coker}(1-C_k)\big)&\mbox{if $\det G=1$}
\end{cases}$ 
\end{enumerate}
\item If $n=1$ and $G$ is an integer dilation matrix, then
\begin{enumerate}
\item $K_0(\cO_{n,G}(\zT^d))=\zZ\oplus\big(\ds\bigoplus_{1\le k\le d,\;\mbox{even}}\mbox{coker}(1-C_k)\big)$
\item $K_1(\cO_{n,G}(\zT^d))=\zZ\oplus\big(\ds\bigoplus_{1\le k\le d,\;\mbox{odd}}\mbox{coker}(1-C_k)\big)$ 
\end{enumerate}
\item If $d=1$, then
\begin{enumerate}
\item $K_0(\cO_{n,m}(\zT))=\zZ_{n-1}$ and $K_1(\cO_{n,m}(\zT))=\zZ_{m-1}$ for $n>1$ and $m\ne0,1$
\item $K_0(\cO_{1,m}(\zT))=\zZ$ and $K_1(\cO_{1,m}(\zT))=\zZ\oplus\zZ_{m-1}$ for $m\ne0,1$
\item $K_0(\cO_{n,1}(\zT))=\zZ\oplus\zZ_{n-1}$ and $K_1(\cO_{n,1}(\zT))=\zZ$ for $n>1$
\end{enumerate} 
\end{enumerate}
\pf For (1), let $d,n>1$ and $\det G\ne1$. Then $C_0=n^d\ne 1$ and $C_d=\det G\ne 1$; 
that is, $1-C_0$ and $1-C_d$ are injective.
Furthermore, for $1\le k\le d-1$, 
$$C_k=n^{d-k}(\det G_{J,I})_{J,I\in\fI_k}.$$ 
Since the characteristic polynomial of a matrix is monic, it follows from Gauss' Lemma that any rational eigenvalue of a matrix in $M_d(\zZ)$ must actually be an integer.
That is,  $\frac{1}{n^{d-k}}\notin\sigma((\det G_{J,I})_{I,J})$ if and only if $1\notin\sigma(C_k).$ 
Thus, $1-C_k$ is injective for each $k=1,..., d-1$ and so, $\ker(1-C_k)=0$ for each $k=0,...,d.$
By Theorem \ref{Alg},
$$K_0(\cO_{n,G}(\zT^d))=\ds\bigoplus_{0\le k\le d,\;\rm even }\mbox{coker}(1-C_k),$$ 
and
$$K_1(\cO_{n,G}(\zT^d))=\ds\bigoplus_{0\le k\le d,\;\rm odd }\mbox{coker}(1-C_k).$$

Now assume $\det G=1$, then, as above, $\ker(1-C_k)=0$ for $k=0,...,d-1.$ But $\ker(1-C_d)=\zZ$ and 
$\mbox{coker}(1-C_d)=\zZ$ whether $d$ is even or odd.
Theorem \ref{Alg} gives
$$K_0(\cO_{n,G}(\zT^d))=\zZ\oplus\big(\ds\bigoplus_{0\le k\le d-1,\;\rm even}\mbox{coker}(1-C_k)\big)$$
and
$$K_1(\cO_{n,G}(\zT^d))=\zZ\oplus\big(\ds\bigoplus_{0\le k\le d-1,\;\rm odd}\mbox{coker}(1-C_k)\big).$$

For (2), let $n=1$ and $\abs{\lambda}>1$ for all eigenvalues $\lambda$ of $G.$ Then $C_0=1$ and $C_d=\det G\ne 1$.
We now wish to show that $\det(1-C_k)\ne0$ for $k=1,...,d$. To this end, choose a basis of $\zC^d$ such that $G$ becomes
upper triangular (not necessarily with integer entries); that is, 
$$G=\begin{pmatrix}
a_{11}&a_{12}&...&a_{1d}\\ 
0&a_{22}&...&a_{2d}\\
\vdots&\ddots&\ddots&\\
0&0&...&a_{dd}
\end{pmatrix}.$$
Then $C_k$ is lower triangular with diagonal entries $\det G_{II}=\prod_{i\in I}a_{ii}\ne1$ (since $\abs{\lambda}>1$
for all $\lambda\in\sigma(G)$) and so $\det(1-C_k)\ne0.$  Hence, $1-C_k$ is injective for $k=1,...,d$ and 
$\ker(1-C_0)=\mbox{coker}(1-C_0)=\zZ$. 
Theorem \ref{Alg} now gives
$$K_0(\cO_{n,G}(\zT^d))=\zZ\oplus\big(\ds\bigoplus_{1\le k\le d,\;\rm even}\mbox{coker}(1-C_k)\big),$$
and
$$K_1(\cO_{n,G}(\zT^d))=\zZ\oplus\big(\ds\bigoplus_{1\le k\le d,\;\rm odd}\mbox{coker}(1-C_k)\big).$$ 

Finally, let us prove (3). If $n>1$ and $m\ne0,1$ then $1-n$ and $1-m$ are injective and 
$\mbox{coker}(1-n)=\zZ_{n-1}:=\zZ/(n-1)\zZ;$ likewise, $\mbox{coker}(1-m)=\zZ_{m-1}$.
Thus, $K_0(\cO_{n,m}(\zT))=\zZ_{n-1}$ and $K_1(\cO_{n,m}(\zT))=\zZ_{m-1}.$ 
If $n=1$ and $m\ne0,1$, then $\ker(1-n)=\mbox{coker}(1-n)=\zZ,$ $\mbox{coker}(1-m)=\zZ_{m-1}$ and $\ker(1-m)=0$   
thus, $K_0(\cO_{1,m}(\zT))=\zZ$ and $K_1(\cO_{1,m}(\zT))=\zZ\oplus\zZ_{m-1}.$
If $n>1$ and $m=1$, then similarly, $K_0(\cO_{n,1}(\zT))=\zZ\oplus\zZ_{n-1}$ and 
$K_1(\cO_{n,1}(\zT))=\zZ.$\\\sq
\end{theorem}

\begin{corollary}\rm Let $d=2$, $n=1$, and $1\notin\sigma(G),$ the spectrum of $G$. Then
\begin{enumerate}
\item if $\det G=1$, then
\begin{enumerate}
\item $K_0(\cO_{1,G}(\zT^2))=\zZ^2$
\item $K_1(\cO_{1,G}(\zT^2))=\zZ^2\oplus\mbox{coker}(1-G^T)$
\end{enumerate}
\item if $\det G\ne 1$, then
\begin{enumerate}
\item $K_0(\cO_{1,G}(\zT^2))=\zZ\oplus\zZ_{1-\det G}$
\item $K_1(\cO_{1,G}(\zT^2))=\zZ\oplus\mbox{coker}(1-G^T)$
\end{enumerate}
\end{enumerate}
\pf Begin with calculating $C_0=n^2=1$, $C_1=nG^T=G^T$ and $C_2=\det G.$ With $1$ not an eigenvalue of $G$, it is guaranteed that $\det(1-G^T)\ne0$ and hence, $1-C_1$ is injective. This leaves us to calculate
\begin{enumerate}
\item $\ker(1-C_0)=\zZ$
\item $\mbox{coker}(1-C_0)=\zZ$
\item $\ker (1-C_1)=0$
\item $\mbox{coker}(1-C_1)=\mbox{coker}(1-G^T)$
\item $\ker(1-C_2)=\begin{cases}\zZ&\mbox{if $\det G=1$}\\0&\mbox{if $\det G\ne1$}\end{cases}$, and
\item $\mbox{coker}(1-C_2)=\begin{cases}\zZ&\mbox{if $\det G=1$}\\\zZ_{1-\det G}&\mbox{if $\det G\ne1$}\end{cases}$
\end{enumerate}\sq
\end{corollary}

We now calculate the K-theory when $F,G\in M_d(\zZ)$ are both diagonal matrices. Let $F=\diag{a_1,...,a_d}$ and $G=\diag{b_1,...,b_d}$ be diagonal integral matrices of non-zero determinant and such that $1\le a_1\le ...\le a_d.$
Let $f$ denote the number of $1$'s in $F;$ that is, $1=a_1=...=a_f< a_{f+1}\le...\le a_d$.
Let $a_I=\prod_{i\in I}a_i$ for $I\in\fI_k$. Then $A_k=\diag{a_I}_{I\in\fI_k}$, 
$B_k=\diag{b_I}_{I\in\fI_k}$ and $$C_k=(\det F)B_kA_k^{-1}=\diag{b_Ia_{I^\prime}}_{I\in\fI_k}.$$
So $\ker(1-C_k)=\zZ^{d_k}$ where $d_k$ is the number of $1$'s in $C_k;$ that is, the number of $I$ making 
$b_Ia_{I^\prime}=1.$ Furthermote, $b_Ia_{I^\prime}=1$ implies $b_I=a_{I^\prime}=1$ and so $\{f+1,...,d\}\subset I.$
So let $v$ be the number of negative ones in $\{b_{f+1},...,b_d\}$ and let $p$ be the number of ones in the same set.
Then
$$d_k=\sum_{r\in\fI} {v\choose 2r}{p\choose (k-d+f)-2r}$$
where $\fI=\{r\in\zN\st 0\le 2r\le v, p-(k-d+f)\le 2r\le k-d+f\}.$

\begin{lemma}\rm With the context of the preceding paragraph, let $p>0$ be the number of ones while $v$ is the number of negative ones. Let $p_k(p,v)$ be the number of combinations of choosing $k$ digits of $1's$ and $-1's$ that multiply to 1 and let $v_k(p,v)$ be similar, but multiplying to $-1$; that is, 
$$p_k(p,v)=\sum_{0,k-p\le 2r\le v,k}{v\choose 2r}{p\choose k-2r}$$
and
$$v_k(p,v)=\sum_{0,k-p\le 2r+1\le v,k}{v\choose 2r+1}{p\choose k-2r-1}$$
for $0<k\le p$, $p_0(p,v)=1$ if $p>0$, and $v_0(p,v)=v_k(p,v)=p_k(p,v)=p_0(0,v)=0$ for $k>p.$ Then
$$\sum_{0\le k\le p} p_k(p,v)=\begin{cases}2^p&\mbox{if $v=0$ and $p>0$}\\ 2^{p+v-1}&\mbox{if $v\ne0$}
\end{cases}.$$
\pf First realize that $p_k(p,v)+v_k(p,v)={p+v\choose k}$ since it is the number of combinations of choosing k digits to multiply to either $1$ or $-1$. Thus, $$\sum_k p_k(p,v)+\sum_kv_k(p,v)=\sum_k{p+v\choose k}=2^{p+v}.$$
So if $v=0$, then the proof is done.

Assume $v\ne0$ and let $v$ be odd. Then claim $v_k(p,v)=p_{p-k}(p,v)$. This can be easily seen by realizing that, for each choice of $k$ digits to multiply to $-1$, the remaining digits multiply to $1$. Thus, $\sum_kv_k(p,v)=\sum_kp_{p-k}(p,v)=\sum_kp_k(p,v)$ and so, $\sum_k p_k(p,v)=\frac{1}{2}2^{p+v}=2^{p+v-1}.$

Now let $v$ be even. Then, for even $k\ne0$
\begin{align*}
v_k(p,v)-v_k(p,v-1)&=\sum_{0,k-p\le 2r+1\le v,k}{v-1\choose 2r+1}{p\choose k-2r-1}\\
&-\sum_{0,k-p+1\le 2r+1\le v,k}{v-1\choose 2r+1}{p\choose k-2r-1}\\
&=\sum_{0,k-p\le 2r+1\le v,k}[{v\choose 2r+1}-{v-1\choose 2r+1}]{p\choose k-2r-1}\\
&=\sum_{0,k-p\le 2r+1\le v,k}{v-1\choose 2r}{p\choose k-2r-1}\\
&=\sum_{0,(k-1)-p\le 2r\le v-1,k-1}{v\choose 2r}{p-1\choose (k-1)-2r}\\
&=p_{k-1}(p,v-1).
\end{align*}
Since $v-1$ is odd and $v_0(p,v)=0$, we get
$$\sum_k v_k(p,v) = \sum_{k>0} v_k(p,v-1)+\sum_{k>0} p_{k-1}(p,v-1)=2^{p+v-2}+2^{p+v-2}=2^{p+v-1}$$
and consequently,
$$\sum_k p_k(p,v)=2^{p+v-1}.$$\sq
\end{lemma}

\begin{theorem}\label{diag}\rm Let $F=\diag{a_1,...,a_d}$ and $G=\diag{b_1,...,b_d}$ be diagonal integral matrices of non-zero determinant and such that $1\le a_1\le ...\le a_d.$ Let $f>0$ denote the number of $1$'s in $F$, let $v$ be the number of negative ones in $\{b_{f+1},...,b_d\}$ and let $p$ be the number of ones in that same set. 
Then
\begin{enumerate}
\item if $p=0$ and $v=0$, then
\begin{enumerate}
\item $K_0(\cO_{F,G}(\zT^d))=\ds\bigoplus_{0\le k\le d,\;\rm even}[\bigoplus_{I\in\fI_k}\zZ_{1-b_Ia_{I^\prime}}]$
\item $K_1(\cO_{F,G}(\zT^d))=\ds\bigoplus_{0\le k\le d,\;\rm odd}[\bigoplus_{I\in\fI_k}\zZ_{1-b_Ia_{I^\prime}}]$ 
\end{enumerate}
\item if $p=0$ and $v>0$, then
\begin{enumerate}
\item $K_0(\cO_{F,G}(\zT^d))=\zZ^{2^{v-1}-1}\oplus\big(\ds\bigoplus_{0\le k\le d,\;\rm odd}[\bigoplus_{I\in\fI_k, b_Ia_{I^\prime}\ne1}\zZ_{1-b_Ia_{I^\prime}}]\big]$
\item $K_1(\cO_{F,G}(\zT^d))=\zZ^{2^{v-1}-1}\oplus\big(\ds\bigoplus_{0\le k\le d,\;\rm odd}[\bigoplus_{I\in\fI_k, b_Ia_{I^\prime}\ne1}\zZ_{1-b_Ia_{I^\prime}}]\big)$ 
\end{enumerate}
\item if $v=0, p>0$, we have
\begin{enumerate}
\item $K_0(\cO_{F,G}(\zT^d))=\zZ^{2^p}\oplus\big(\ds\bigoplus_{0\le k\le d,\;\rm even}[\bigoplus_{I\in\fI_k, b_Ia_{I^\prime}\ne1}\zZ_{1-b_Ia_{I^\prime}}]\big)$
\item $K_1(\cO_{F,G}(\zT^d))=\zZ^{2^p}\oplus\big(\ds\bigoplus_{0\le k\le d,\;\rm odd}[\bigoplus_{I\in\fI_k, b_Ia_{I^\prime}\ne1}\zZ_{1-b_Ia_{I^\prime}}]\big)$ 
\end{enumerate}
\item if  $v\ne0, p>0$, then
\begin{enumerate}
\item $K_0(\cO_{F,G}(\zT^d))=\zZ^{2^{p+v-1}}\oplus\big(\ds\bigoplus_{0\le k\le d,\;\rm even}[\bigoplus_{I\in\fI_k, b_Ia_{I^\prime}\ne1}\zZ_{1-b_Ia_{I^\prime}}]\big)$
\item $K_1(\cO_{F,G}(\zT^d))=\zZ^{2^{p+v-1}}\oplus\big(\ds\bigoplus_{0\le k\le d,\;\rm odd}[\bigoplus_{I\in\fI_k, b_Ia_{I^\prime}\ne1}\zZ_{1-b_Ia_{I^\prime}}]\big)$ 
\end{enumerate}
\end{enumerate}
\pf  Begin by calculating
\begin{enumerate}
\item $\ds\bigoplus_{k,\,\rm even (odd)}\ker (1-C_k)=\bigoplus_{k,\,\rm even (odd)}\zZ^{d_k}$
\item $\ds\bigoplus_{k,\,\rm even (odd)}\mbox{coker}(1-C_k)=\bigoplus_{k,\,\rm even (odd)}\big(\ds\bigoplus_{b_Ia_{I^\prime}\ne1, I\in\fI_k} \zZ_{1-b_Ia_{I^\prime}}\big)\oplus\big(\bigoplus_{k,\,\rm even (odd)}\zZ^{d_k}\big).$
\end{enumerate}
Thus,
\begin{align*}
K_0(\cO_{F,G}(\zT^d))&=\big(\bigoplus_{0\le k\le d, \,\rm even} \mbox{coker}(1-C_k)\big)
\oplus\big(\bigoplus_{0\le k\le d, \,\rm odd} \ker(1-C_k)\big)\\
&=\big(\bigoplus_{k,\,\rm odd}\zZ^{d_k}\big)\oplus\big(\ds\bigoplus_{k,\,\rm even}\big(\ds\bigoplus_{b_Ia_{I^\prime}\ne1, I\in\fI_k} \zZ_{1-b_Ia_{I^\prime}}\big)\oplus\big(\bigoplus_{k,\,\rm even}\zZ^{d_k}\big)\\
&=\big(\bigoplus_k\zZ^{d_k}\big)\oplus\big(\ds\bigoplus_{k,\,\rm even}\big(\ds\bigoplus_{b_Ia_{I^\prime}\ne1, I\in\fI_k} \zZ_{1-b_Ia_{I^\prime}}\big)\\
&=\zZ^{\sum_k d_k}\oplus\big(\ds\bigoplus_{k,\,\rm even}\big(\ds\bigoplus_{b_Ia_{I^\prime}\ne1, I\in\fI_k} \zZ_{1-b_Ia_{I^\prime}}\big).\\
\end{align*}
Similarly for $K_1;$
$$K_1(\cO_{F,G}(\zT^d)=\zZ^{\sum_k d_k}\oplus\big(\ds\bigoplus_{k,\,\rm odd}\big(\ds\bigoplus_{b_Ia_{I^\prime}\ne1, I\in\fI_k} \zZ_{1-b_Ia_{I^\prime}}\big).$$
If $p=0$ and $v=0$, then there is no way to multiply to get $1$. Hence, $1-C_k$ is injective for all $k$ and the result follows readily. If $p=0$ and $v>0$, then $d_k=0$ for all odd $k$ or $k> v$ and $d_k={v\choose k}$ for all even $k\le v.$
Thus, 
$$\sum_k d_k=\sum_{2\le 2k\le v}{v\choose 2k}=(\sum_{0\le 2k\le v}{v\choose 2k})-1=2^{v-1}-1$$
and the rest follows.

If $p>0$, then 
$$d_k=\sum_{r\in\fI} {v\choose 2r}{p\choose (k-d+f)-2r}$$
where $\fI=\{r\in\zN\st 0\le 2r\le v, p-(k-d+f)\le 2r\le k-d+f\}.$ Let $k^\prime=k-d+f$ and then
$$d_k=p_{k^\prime}(p,v)$$
where $p_{k^\prime}(p,v)$ was defined in the above lemma and so, 
$$\sum_k d_k=\sum_k p_k(p,v)=\begin{cases} 2^p&\mbox{if $v=0$}\\ 2^{p+v-1}&\mbox{if $v\ne0$}
\end{cases}$$
and once again the rest follows.\\\sq
\end{theorem}

\begin{remark}\rm If $f,$ the number of ones in $F,$ is $0$ then $C_k$ is injective for $k=0,...,d-1$. It is then very simple to calculate the $K$-groups whether $\det G=1$ or $\det G\ne1$.
\end{remark}

\begin{corollary}\rm If $F=n1_d, G=m1_d\in M_d(\zZ)$ where $n\in\zN$ and $m\in\zZ$ are non-zero, then 
\begin{enumerate}\item if either $n>1$ or $\abs{m}\ne1$, then
\begin{enumerate}
\item $K_0(\cO_{n,m}(\zT^d))=\ds\bigoplus_{0\le k\le d,\;\rm even}\zZ_{1-n^{d-k}m^k}^{d\choose k}$
\item $K_1(\cO_{n,m}(\zT^d))=\ds\bigoplus_{0\le k\le d,\;\rm odd}\zZ_{1-n^{d-k}m^k}^{d\choose k}$ 
\end{enumerate}
\item if $n=m=1$, then $K_0(\cO_{1,1}(\zT^d))=K_1(\cO_{1,1}(\zT^d))=\zZ^{2^d}$ 
\item if $n=1$, $m=-1$, then $K_0(\cO_{1,-1}(\zT^d))=K_1(\cO_{1,-1}(\zT^d))=\zZ^{2^{d-1}}\oplus \zZ_{2}^{2^{d-1}}$ 
\end{enumerate}
\end{corollary}

\section{Acknowledgements}This paper is the second product of my doctoral dissertation at the University of Calgary. It
would not have been possible without Dr. Berndt Brenken for his 
insight on corrections and enhancements of the material discussed here.

I am also indebted to NSERC, the Department of Mathematics and Statistics at the
University of Calgary and the Department of Mathematics and Statistics at the University of Regina 
in providing funding to finance my mathematical studies.


\end{document}